\documentclass{amsart}
\usepackage{graphicx,amssymb}
\usepackage{mathtools}
\usepackage{amsmath,  amsfonts, latexsym}
\usepackage{verbatim}

\usepackage{color}

\newtheorem{theorem}{Theorem}[section]

\theoremstyle{definition}

\theoremstyle{remark}

\numberwithin{equation}{section}

\def\XXint#1#2#3{{\setbox0=\hbox{$#1{#2#3}{\int}$ }
\vcenter{\hbox{$#2#3$ }}\kern-.6\wd0}}

\newcommand{\divx}{\operatorname{div}}
\newcommand{\nablax}{\nabla}
\newcommand{\tr}{\mathrm{tr}\,}
\newcommand{\va}{\vec{a}}
\newcommand{\vb}{\vec{b}}
\newcommand{\vv}{\vec{v}}
\newcommand{\nn}{\vec{n}}
\newcommand{\vw}{\vec{w}}

\newcommand{\AAA}{\mathbb{A}}
\newcommand{\BB}{\mathbb{B}}
\newcommand{\TS}{\mathbb{S}}
\newcommand{\BBkpt}{\mathbb{B}_{k_p(t)}}

\newcommand{\DD}{\mathbb{D}}
\newcommand{\CCkpt}{\mathbb{C}_{k_p(t)}}
\newcommand{\DDkpt}{\mathbb{D}_{k_p(t)}}
\newcommand{\LLkpt}{\mathbb{L}_{k_p(t)}}

\newcommand{\FFkpt}{\mathbb{F}_{k_p(t)}}
\newcommand{\GG}{\mathbb{G}}
\newcommand{\II}{\mathbb{I}}

\newcommand{\OO}{\mathbb{O}}
\newcommand{\LL}{\mathbb{L}}

\newcommand{\TT}{\mathbb{T}}
\newcommand{\Sstar}{\overset{*}{\mathbb{S}}}
\newcommand{\Dstar}{\overset{*}{\mathbb{D}}}
\newcommand{\Astar}{\overset{*}{\mathbb{A}}}
\newcommand{\TTel}{\TT_{\mathrm{el}}}

\newcommand{\bS}{\mathbb{S}}
\newcommand{\0}{\mathbb{O}}

\newcommand{\Aold}{\overset{\triangledown}{\mathbb{A}}}
\newcommand{\Acirc}{\overset{\circ}{\mathbb{A}}}
\newcommand{\Asq}{\overset{\square}{\mathbb{A}}}

\newcommand{\pthns}{p_{\mathrm{th}}^{\mathrm{NS}}}

\newcommand{\qe}{\vec{j}_e}
\newcommand{\qn}{\vec{j}_{\eta}}

\newcommand{\dd}{\,{\rm d}}

\begin{document}

\title[Viscoelastic rate-type fluids with stress-diffusion]{PDE analysis of a class of
thermodynamically compatible viscoelastic rate-type fluids with stress-diffusion} \thanks{The authors acknowledge the support of the ERC-CZ project LL1202, financed by M\v{S}MT}


\author[M.~Bul\'{\i}\v{c}ek]{Miroslav Bul\'{\i}\v{c}ek}
\address{Mathematical Institute, Charles University, Faculty of Mathematics and Physics, Sokolovsk\'a 83, 186 75 Prague 8, Czech Republic}
\email{mbul8060@karlin.mff.cuni.cz}
\thanks{}

\author[J.~M\'{a}lek]{Josef M\'{a}lek}
\address{Mathematical Institute, Charles University, Faculty of Mathematics and Physics, Sokolovsk\'a 83, 186 75 Prague 8, Czech Republic}
\email{malek@karlin.mff.cuni.cz}
\thanks{}

\author[V.~Pr\accent23u\v{s}a]{V\'{\i}t Pr\accent23u\v{s}a}
\address{Mathematical Institute, Charles University, Faculty of Mathematics and Physics,  Sokolovsk\'a 83, 186 75 Prague 8, Czech Republic}
\email{prusv@karlin.mff.cuni.cz}
\thanks{}

\author[E.~S\"{u}li]{Endre S\"{u}li}
\address{Mathematical Institute, University of Oxford, Woodstock Road, Oxford OX2 6GG, United Kingdom}
\email{Endre.Suli@maths.ox.ac.uk}
\thanks{}

\subjclass[2000]{Primary 35Q35, 76A05, 76A10}

\date{}

\begin{abstract}
We establish the long-time existence of large-data weak solutions to a system of nonlinear partial differential equations. The system of interest governs the motion of non-Newtonian fluids described by a simplified viscoelastic rate-type model with a stress-diffusion term. The simplified model shares many qualitative features with more complex viscoelastic rate-type models that are frequently used in the modeling of fluids with complicated  microstructure. As such, the simplified model provides important preliminary insight into the mathematical properties of these more complex and practically relevant models of non-Newtonian fluids. The simplified model that is analyzed from the mathematical perspective is shown to be thermodynamically consistent, and we extensively comment on the interplay between the thermodynamical background of the model and the mathematical analysis of the corresponding initial-boundary-value problem.
\end{abstract}

\maketitle

\section{Introduction}

The main goal of this study is to establish the long-time existence of large-data weak solutions to the following simplified set of governing equations encountered in the mechanics of incompressible non-Newtonian fluids. This simple model shares many standard properties with more complex viscoelastic rate-type fluid models with stress-diffusion, used in applications.

For any given Lipschitz domain $\Omega \subset \mathbb{R}^d$, $d\ge 2$, $\vv_0: \Omega \to \mathbb{R}^d$, $b_0:\Omega\to \mathbb{R}$ and for any  $T>0$, we set $Q:=\Omega \times (0,T)$ and we seek the functions $(\vv,p,b): Q \to \mathbb{R}^d\times\mathbb{R}\times\mathbb{R}$ satisfying, in $Q$, the system of PDEs:
\begin{align}
\divx \vv&=0, \label{eq10}\\
\varrho(\partial_t \vv + \divx(\vv \otimes \vv))- \divx \TT&=0, \label{eq10.5} \\
\TT &=-p\,\II + 2\nu\DD - \sigma (\nablax b \otimes \nablax b), \label{eq11} \\
\nu_1 \partial_t b + \nu_1 \divx (b\vv) + \mu (b^2 - b) - 2\sigma b^2 \Delta b  &= 0, \label{eq12}
\end{align}
together with the following boundary conditions, where $\nn$ denotes the unit outward normal vector on $\partial \Omega$:
\begin{equation}\label{eq.bc}
\vv=\vec{0} \quad \textrm{ and } \quad \nabla b \cdot \nn =0 \qquad \textrm{ on } \partial \Omega \times (0,T),
\end{equation}
and the initial conditions:
\begin{equation}
\vv(0,\cdot)=\vv_0(\cdot) \quad \textrm{ and } \quad b(0,\cdot)=b_0(\cdot) \qquad \textrm{ in } \Omega.
\label{eq.ic}
\end{equation}
Here, $\vv$ is the velocity, $\otimes$ stands for the dyadic product ($(\va\otimes\vb)_{ij} = a_i b_j$),  $\TT$ is the Cauchy stress, $p$ is the spherical stress (the modified pressure), $b$ is a scalar quantity characterizing the volumetric elastic changes exhibited by the fluid;
$\DD$ denotes the symmetric part of the velocity gradient, i.e., $\DD=(\nabla \vv + (\nabla \vv)^{\textrm{T}})/2$; $\nu$ is the viscosity, $\varrho$ is the density, and $\nu_1$, $\sigma$ and $\mu$ are other material constants, all of which are positive.

Provided further that (the definitions of all function spaces appearing in the paper are given in Section \ref{FS})
\begin{equation}
\begin{split}
\vv_0 &\in L^2_{0,\divx}(\Omega)^d, \\
b_0 &\in W^{1,2}(\Omega), \quad  b_0>0 ~\textrm{ a.e. in } \Omega, \quad b_0 \in L^{\infty}(\Omega) \quad \textrm{ and } \quad \frac{1}{b_0}\in L^{\infty}(\Omega)\,,
\end{split} \label{data1}
\end{equation}
we will prove, see Section \ref{Proof}, the following statement:

\medskip

\begin{minipage}{0.85\hsize}
\emph{For any Lipschitz domain $\Omega \subset \mathbb{R}^d$, $T>0$, $\nu>0$, $\nu_1>0$, $\mu>0$, $\sigma>0$, $\varrho>0$ and  for any $(\vv_0, b_0)$ satisfying \eqref{data1}, there exists a couple $(\vv,b)$ solving the problem \eqref{eq10}--\eqref{eq.ic} in a weak sense.}
\end{minipage}

\bigskip

\noindent The precise definition of \textit{weak solution} is given in Section \ref{FS}.

Note that if $\sigma=0$ in the above equations the problem splits into two separate problems: a transport equation with damping for $b$ and the standard Navier--Stokes equations for $(\vv,p)$. Since the work of Leray \cite{leray.j:sur}  on the incompressible Navier--Stokes equations, see also \cite{hopf51, caffarelli.l.kohn.r.ea:partial}, the question of long-time existence of large-data weak solutions has also been explored for more general classes of viscous fluids; further references will be given in the next section. It is then natural to attempt to advance this program by exploring how large the class of fluids might be for which long-time and large-data existence of weak solutions can be established. This task is particularly interesting if one considers fluid models that include terms that are associated with elastic properties of the material. This paper aims to contribute to the study of this question.

Note also that if $\sigma>0$ then not only is there a diffusion term present in the equation for $b$, but also the Korteweg stress appears in the expression for the Cauchy stress $\TT$ featuring in the equation for the balance of linear momentum. This structure of the governing equations results from a careful derivation of the model based on the thermodynamical approach established for viscoelastic rate-type fluids in \cite{Rajagopal2000} and further refined in \cite{malek.j.rajagopal.kr.ea:on} (see also \cite{malek.j.prusa.v:derivation}), and extended to rate-type fluids with stress-diffusion in \cite{MPSS17}. This thermodynamical approach automatically guarantees that the resulting model is consistent with the laws of thermodynamics. From the point of view of PDE analysis, the approach readily provides the relevant \textit{a priori} estimates upon which the analysis is based, and which would be otherwise completely nontrivial to discover solely from the PDE system \eqref{eq10}--\eqref{eq12}.

Thus, the second goal of this paper is to highlight how the analysis of PDEs for complex systems such as \eqref{eq10}--\eqref{eq12} is related to the thermodynamical approach used in the derivation of the model and to show that the problem \eqref{eq10}--\eqref{eq.ic} therefore merits investigation.

The paper is structured as follows. In the next section, we place the model under consideration into a hierarchy of phenomenological models. In Section \ref{der}, we sketch its derivation. As was noted above, this thermodynamical approach automatically guarantees the consistency with the second law of thermodynamics and provides directly the \textit{a priori} estimates available for the system.  In Section \ref{FS}, we re-derive these \textit{a priori} estimates, but now from a purely PDE analytical point of view, and then, after introducing the appropriate function spaces, we precisely state the main result of the paper. The proof of the theorem is contained in Section \ref{Proof}. The existence of solutions to the finite-dimensional approximating system, upon which the proof of the main result is based, is given in the Appendix.

\section{Incompressible non-Newtonian fluids: a brief overview}


For incompressible fluids, the Cauchy stress tensor $\TT$ decomposes as
\begin{align}
\TT = - \phi \II + \TS, \label{j1}
\end{align}
where $\phi$ is a scalar unknown quantity and $\TS$ is related to the symmetric part of the velocity gradient $\DD$ via the so-called constitutive relation. Three classes of constitutive relations are frequently used; they are of the following forms:
\begin{align}
\GG(\TS,\DD) &= \0, \label{j2} \\
\GG(\Sstar,\TS,\Dstar,\DD) &= \0, \label{j3} \\
\GG(\Sstar,\TS,\Dstar,\DD) - \Delta \TS &= \0. \label{j4}
\end{align}
Here $\GG$ stands for an arbitrary continuous tensorial function and $\Astar$ signifies an objective time derivative.

We will give examples and discuss briefly the usefulness of these classes. Then we shall look at the associated initial-boundary-value problems for internal flows from the perspective of the long-time existence of large-data weak solutions.

Besides the (linear, i.e. Newtonian) Navier--Stokes fluid for which the constitutive relation is given by $\TS=2\nu \DD$ with $\nu>0$, the class \eqref{j2} includes for example the Bingham fluid, described by $2\nu\DD=\frac{(|\TS|- \tau_*)^{+}}{|\TS|} \TS$, as well as various generalizations of activated or non-activated (stress) power-law fluids that can be characterized by
$$
2\nu(|\TS|^2, |\DD|^2)(|\DD|- d_*)^{+}\DD=2\alpha(|\TS|^2, |\DD|^2) (|\TS|- \tau_*)^{+}\TS,
$$
which involves two non-negative activation parameters $d_*$ and $\tau_*$, but only one is assumed to be positive in a given model. The possibility of \eqref{j2} to incorporate nonlinear relations between $\TS$ and $\DD$ and to include sudden changes from one type of response to another one (within the considered class) makes the class \eqref{j2} suitable for modeling even mixing type phenomena, 
and the class \eqref{j2} is therefore very popular in many areas (see \cite{MRR95, bulcek.m.gwiazda.p.ea:on*3} and \cite{malek.j.prusa.v:derivation} for details and further references).

On the other hand, fluids described by the class \eqref{j2} are not capable of capturing fundamental phenomena such as stress relaxation or nonlinear creep observed in most fluids with complex microstructure. If these phenomena are of interest, then the class \eqref{j3} is a suitable choice.

As was already mentioned above in connection with \eqref{j3}, $\overset{*}{\mathbb{A}}$ denotes an objective time derivative. The set of objective derivatives includes the upper-convected Oldroyd derivative defined by
\begin{equation*}
\Aold  := \frac{\dd \AAA}{\dd t} - \LL\AAA-\AAA\LL^{\textrm{T}} \quad \textrm{ where } \quad \LL:=\nabla \vv,
\end{equation*}
the Jaumann--Zaremba (corotational) derivative specified by
\begin{equation*}
\Acirc := \frac{\dd \AAA}{\dd t} - \mathbb{W}\AAA-\AAA\mathbb{W}^{\textrm{T}} \quad \textrm{ where } \quad \mathbb{W}:=(\LL - \LL^{\textrm{T}})/2,
\end{equation*}
or the Gordon--Schowalter derivative defined by
\begin{equation*}
\Asq  := \Acirc - a(\DD\AAA-\AAA\DD) \quad \textrm{ where } a\in[-1,1].
\end{equation*}
Assuming that $\TT = - \phi \II + 2\nu\DD + \AAA$, we see that the appropriate choice of material parameter $\nu$ and the objective derivative lead to the popular Maxwell, Oldroyd-B and Johnson--Segalman models respectively:
\begin{equation*}
\begin{split}
\tau\Aold + \AAA & = 2\nu_1\DD \qquad \textrm{ with } \quad \nu=0 \textrm{ and } \tau=\frac{\nu_1}{E}, \\
\tau\Aold + \AAA & = 2\nu_1\DD \qquad \textrm{ with } \quad \nu>0 \textrm{ and } \tau=\frac{\nu_1}{E}, \\
\tau\Asq+\TS & = 2a\DD \qquad \,\, \textrm{ with } \quad \nu>0 \textrm{ and } a\in[-1,1].
\end{split}
\end{equation*}
Although these models are widely used, there are several subtle issues regarding their physical underpinnings. These include the ambiguity with respect to the choice of the objective derivative, the possibility of the derivation of the model at a
purely macroscopic level, the consistency of the models with the second law of thermodynamics, the extension of the models to the compressible setting, and the inclusion of thermal effects into the models.

To address some of these issues Rajagopal and Srinivasa \cite{Rajagopal2000} provided a simple, yet general method for the derivation of thermodynamically consistent rate-type fluid models. The method is based on the concept of \textit{natural configuration} and on the knowledge of constitutive relations for \emph{two scalar quantities}: the Helmholtz free energy (characterizing how the material stores energy) and the rate of entropy production (characterizing how the material dissipates energy).  Using the method, it is possible to derive new thermodynamically compatible classes of \emph{non-linear} viscoelastic rate-type fluid models, and to investigate the conditions under which these models reduce to the standard Oldroyd-B, Maxwell and Burgers models. More recently, the approach has been refined and extended to compressible fluids and thermal processes (see \cite{malek.j.rajagopal.kr.ea:on, malek.j.prusa.v:derivation, hron.j.milos.v.ea:on}).

Rate-type fluid models with stress-diffusion, c.f. \eqref{j4}, are popular in the modelling of shear and vorticity banding phenomena; see, for example, the reviews \cite{dhont.jg.briels.w:gradient, fardin.ma.ober.tj.ea:potential, divoux.t.fardin.ma.ea:shear}.  The \textit{ad hoc} addition of the regularizing term $-\sigma \Delta \TS$ to the rate-type equation \eqref{j3} raises doubts about the consistency of the resulting model with the second law of thermodynamics and calls for the specification of appropriate boundary conditions for $\TS$. These issues are addressed in the recent study \cite{MPSS17}. Here, in Section \ref{der}, we sketch the derivation under several simplifying assumptions.

As to the question of existence of solutions, for any time interval and any reasonable set of data, a satisfactory theory, within the context of weak solutions, has been established during the last fifty years for a general class of fluids given by \eqref{j2} with a continuous monotone function $\GG$ (see \cite{leray.j:sur, caffarelli.l.kohn.r.ea:partial, kiselev-ladyzh, ladyzhenskaya.oa:on, ladyzhenskaya.oa:mathematical*1, lions.jl:quelques, bulcek.m.ettwein.f.ea:on,
BBN93, MNR93, MNR00, FMS00, bulcek.m.malek.j.ea:naviers, Wolf07, DRW10, BDS13,  bulcek.m.gwiazda.p.ea:on*3, BM16, MaZa17}).

Regarding the long-time and large-data existence theory for rate-type fluids in three dimensions
much less is known (see \cite{lions.pl.masmoudi.n:global, hu.d.lelievre.t:new, masmoudi.n:global} or \cite{LeBrisLelievre}).

One expects that the addition of the term $-\Delta \TS$ to a rate-type equation will automatically improve its mathematical properties and a well-founded theory for a class of fluids of the type \eqref{j4} will therefore emerge. However, most of the results that are in place concern two-dimensional (or steady) flows, see \cite{el-kareh.aw.leal.lg:existence, barrett.jw.boyaval.s:existence, constantin.p.kliegl.m:note, chupin.l.martin.s:stationary, LMN15, ER2015} or require stronger regularization, see \cite{kreml-pokorny-salom}.

Unlike the above two-dimensional results, the general analytical approach developed and presented in Sections \ref{FS} and \ref{Proof} concerns $d$-dimensional flows. Although the tensor $\TS$ is ultimately ``replaced" here by a scalar quantity $b$, the analysis is performed in such a way that it avoids properties which are unlikely to be true for a general tensor $\TS$.

\section{Derivation of the model} \label{der}

We denote by the superscript $\dot{~}$ the material derivative, i.e., $\dot{z} = \tfrac{\partial z}{\partial t} + (\vv\cdot \nablax) z$ for a scalar function $z$; if $z$ stands for a longer expression, we write $\dot{\overline{z}}$ instead of $\dot{z}$ in order to avoid confusion. For vector- or tensor-valued quantities, the same notation is used and the above formula is applied to each component.

The fundamental balance equations of mass, linear momentum and energy as well as the formulation of the second law of thermodynamics are the following:
\begin{align}
\dot{\varrho} &= - \varrho\, \divx \vv,\label{p1}\\
\varrho \dot{\vv} &= \divx \TT, \qquad \TT = \TT^{\textrm{T}},\label{p2}\\
\varrho \dot{e} &= \TT : \DD - \divx \qe, \label{p3}\\
\varrho \dot{\eta} &= \varrho \zeta - \divx \qn \quad \textrm{ with } \zeta \ge 0,\label{p4}
\end{align}
where $\varrho$ is the density, $\vv$ is the velocity, $e$ is the specific internal energy, $\eta$ is the specific entropy, $\DD= (\nablax \vv + (\nabla \vv)^{\textrm{T}})/2$ is the symmetric part of the velocity gradient, $\TT$ is the Cauchy stress tensor, $\qe$ and $\qn$ are the energy and entropy fluxes, and $\zeta$ stands for the specific rate of entropy production.

Introducing the Helmholtz free energy $\psi$ by
\begin{equation*}
\psi:= e - \theta \eta,
\end{equation*}
where $\theta$ stands for the (positive) temperature, the equations \eqref{p3} and \eqref{p4} lead to
\begin{equation*}
\TT : \DD - \varrho \dot{\psi} - \divx (\qe - \theta \qn) - \varrho \eta\dot{\theta} - \nablax \theta \cdot \qn = \varrho \theta \zeta
\quad \textrm{ with } \zeta \ge 0.
\end{equation*}
In what follows, we restrict ourselves to isothermal processes (referring the reader to \cite{MPSS17} for the development of models where both mechanical and thermal processes are included). Consequently, $\theta>0$ is constant and the last identity reduces to
\begin{equation}
\TT : \DD - \varrho \dot{\psi} - \divx (\qe - \theta \qn)  = \xi \quad \textrm{ with } \xi \ge 0, \label{p5}
\end{equation}
where $\xi:=\varrho \theta \zeta$ denotes the rate of dissipation. Note that \eqref{p5} simplifies further to $\TT : \DD - \varrho \dot{\psi} = \xi$ provided that $\qn= \tfrac{\qe}{\theta}$, which we do not require here, however.

The approach that we will exploit is based on the concept of natural configuration that splits the total response described by the deformation tensor $\mathbb{F}$ between the current and initial configuration into the purely elastic (reversible, non-dissipative) part $\FFkpt$ that operates between the natural and current configuration and the dissipative part $\mathbb{G}$ that maps from the reference to the natural configuration, i.e. $\mathbb{F} = \FFkpt \mathbb{G}$. In analogy with the relations $\LL= \dot{\mathbb{F}}\mathbb{F}^{-1}$ and $\DD= (\LL+ \LL^{\textrm{T}})/2$, we set $\LLkpt:= \dot{\mathbb{G}} \mathbb{G}^{-1}$ and define $\DDkpt := (\LLkpt + (\LLkpt)^{\textrm{T}})/2$. We also set $\BBkpt:=\FFkpt \FFkpt^{{\textrm{T}}}$.

In order to derive the constitutive relations for the Cauchy stress tensor and the energy and entropy fluxes we start with postulating the constitutive relation for the Helmholtz free energy in the form:
\begin{equation}
\psi = \psi_0(\varrho) + \frac{\mu}{2\varrho}(\tr \BBkpt - 3 - \ln \det \BBkpt) + \frac{\sigma}{2\varrho}|\nablax\, \tr \BBkpt|^2 =: \psi_0(\varrho) + \frac{\pi}{\varrho}, \label{5a}
\end{equation}
with $\mu>0$, $\sigma>0$ constant, and
\[\pi: = \frac{\mu}{2}(\tr \BBkpt - 3 - \ln \det \BBkpt) + \frac{\sigma}{2}|\nablax\, \tr \BBkpt|^2.  \]
The following identities (see \cite{malek.j.rajagopal.kr.ea:on} for their derivation) hold:
\begin{align}
\dot{\overline{\BBkpt}} &= \LL \BBkpt + \BBkpt \LL^{\textrm{T}} - 2\, \FFkpt \DDkpt \FFkpt^{\textrm{T}}, \label{p6}\\
\dot{\overline{\tr \BBkpt}} &= 2\, \BBkpt : \DD - 2\, \CCkpt : \DDkpt,\label{p7}\\
\dot{\overline{\ln \det \BBkpt}} & = 2\, \II : \DD - 2\,\II: \DDkpt, \label{p8}\\
\frac{1}{2}\dot{\overline{|\nablax \tr \BBkpt|^2}} &=  \nablax \tr \BBkpt \cdot \left(\nablax \dot{\overline{\tr \BBkpt}} - (\nablax \vv)\, \nablax \tr \BBkpt\right). \label{p9}
\end{align}
Setting
\[
b:= \tr \BBkpt
\]
and using the formula
\[
  \nablax b \cdot \nablax \dot{b} = \divx (\dot{b}  \, \nablax b ) - \Delta b \, \dot{b},
\]
the above identities \eqref{5a}, \eqref{p7} and \eqref{p9} take the form
\begin{align}
\psi &= \psi_0(\varrho) + \frac{\pi}{\varrho} \quad \textrm{ with } \quad \pi = \frac{\mu}{2}(b-3 - \ln \det \BBkpt) + \frac{\sigma}{2}|\nablax\, b|^2, \label{p10} \\
\dot{b} &= 2\, \BBkpt : \DD - 2\, \CCkpt : \DDkpt,\label{p11} \\
\frac{\dot{\overline{|\nablax b|^2}}}{2} &=  -\left( (\nablax b \otimes \nablax b) + 2 \Delta b  \,  \BBkpt\right) :\DD + \divx(\dot{b}\,\nablax b) + 2 \Delta b  \,  \CCkpt : \DDkpt. \label{p12}
\end{align}
Consequently, using also \eqref{p1} and \eqref{p8}, we compute
\begin{align*}
\varrho \dot{\psi} &= \varrho\psi_0'(\varrho)\dot{\varrho} + \dot{\pi} - \frac{\pi \dot{\varrho}}{\varrho} \\
&= \left( \mu(\BBkpt - \II) - \sigma ( (\nablax b \otimes \nablax b) + 2 \Delta b  \,  \BBkpt )\right) :\DD \\
&\quad - \left(\mu (\CCkpt - \II) - 2\sigma \Delta b  \,  \CCkpt \right) : \DDkpt + \divx(\sigma \dot{b}\,\nablax b) - (\varrho^2 \psi_0'(\varrho) - \pi)\divx\vv.
\end{align*}
Inserting the result into \eqref{p5} we obtain
\begin{align*}
\xi&=\left(\TT - \mu(\BBkpt - \II) + \sigma ( (\nablax b \otimes \nablax b) + 2 \Delta b  \,  \BBkpt ) + (\varrho^2 \psi_0'(\varrho) - \pi)\II \right) :\DD \\
&\quad + \left(\mu (\CCkpt - \II) - 2\sigma \Delta b  \,  \CCkpt \right) : \DDkpt - \divx(\qe - \theta \qn + \sigma \dot{b}  \, \nablax b) \quad \textrm{ with } \quad \xi \ge 0.
\end{align*}
Finally, setting\footnote{Note that the first line states how the entropy flux $\qn$ is related to the energy flux $\qe$, while in the second and third lines we have introduced new notation to simplify the formulae.}
\begin{align*}
\qn &= \frac{\qe +  \sigma \dot{b}  \, \nablax b}{\theta}, \\
\pthns &:= \varrho^2\, \psi_0'(\varrho), \\
\TTel &:= - (\pthns - \pi)\II + \mu(\BBkpt - \II) - \sigma ( (\nablax b \otimes \nablax b) + 2 \Delta b  \,  \BBkpt ),
\end{align*}
we arrive at
\begin{align}
\xi&=\left(\TT - \TTel \right) :\DD + \left(\mu (\CCkpt - \II) - 2\sigma \Delta b  \,  \CCkpt \right) : \DDkpt \quad \textrm{ with } \quad \xi \ge 0, \label{p13}
\end{align}
which provides key information in the derivation of an entire hierarchy of models, ranging from (compressible/incompressible) Euler and Navier--Stokes fluids, through Maxwell and Oldroyd-B fluids, to diffusive Maxwell and diffusive Oldroyd-B fluids. In addition, the procedure also clarifies how these fluids store and dissipate energy and consequently what quantities are \textit{a priori} controlled by the initial data (which then indicates what \textit{a priori} information/estimates one can use in the analysis of the model).

To illustrate the procedure and to derive a simplified (toy) model that is then analyzed in the second part of this paper, we rewrite \eqref{p13} in the following way. Using the subscript $_\delta$ to signify the deviatoric (traceless) part of a tensor, i.e. $\AAA_\delta:= \AAA - \tfrac13\tr \AAA\, \II$, noticing that
\[ \AAA : \BB = \AAA_\delta:\BB_\delta + \tr\AAA \frac{\tr\BB}{3} \]
and denoting the \textit{mean normal stress} by $m$, i.e.,
\[ m: = \frac{1}{3} \tr \TT,\]
we rewrite \eqref{p13} as
\begin{equation}
\begin{split}
\xi=\left(\TT_\delta - [\TTel]_\delta \right) :\DD_\delta &+ \bigg(m - \frac{\tr \TTel}{3}\bigg) \divx \vv \\
&+\left(\mu [\CCkpt]_\delta - 2\sigma \Delta b  \,  [\CCkpt]_\delta \right) : [\DDkpt]_\delta \\
&+ \left(\mu (b - 3) - 2\sigma b\,\Delta b \right)\frac{\tr \DDkpt}{3} \quad \textrm{ with } \quad \xi \ge 0,
\end{split}\label{p14}
\end{equation}
where
\begin{align*}
[\TTel]_\delta &= \mu[\BBkpt]_\delta - \sigma ( (\nablax b \otimes \nablax b)_\delta + 2 \Delta b  \,  [\BBkpt]_\delta ), \\
\frac{\tr \TTel}{3} &= - \pthns + \frac13 [\mu (b-3) - \sigma |\nablax\, b|^2 - 2\sigma b\,\Delta b] + \pi \\
&\qquad \qquad \qquad \textrm{ with } \pi = \frac{\mu}{2}(b-3 - \ln \det \BBkpt) + \frac{\sigma}{2}|\nablax\, b|^2.
\end{align*}
Referring to \cite{malek.j.prusa.v:derivation} for details, we wish to emphasize that the structure of \eqref{p14} is rich enough to incorporate the compressible and incompressible Euler and Navier--Stokes, Maxwell, Oldroyd-B, Giesekus,  diffusive Maxwell, diffusive Oldroyd-B, and diffusive Giesekus models.

Here, we make three simplifications: the fluid is assumed to be incompressible, the density is taken to be uniform and the elastic part of the deformation is supposed to be purely spherical\footnote{This is reminiscent of the response leading to the Euler fluid.}. This means that
\begin{equation}
\divx \vv = 0, \quad \varrho \textrm{ is a positive constant} \quad \textrm{ and } \quad [\CCkpt]_{\delta} = \OO. \label{ass1}
\end{equation}
Then,
\[\CCkpt = [\CCkpt]_\delta + \frac{\tr \CCkpt}{3} \II = \OO + \frac{\tr \BBkpt}{3}\, \II = \frac{b}{3}\, \II,\]
where $\OO$ is the zero tensor; note that $[\BBkpt]_\delta$ also vanishes. As a consequence of these simplifications, \eqref{p11} and \eqref{p14} reduce to
\begin{align}
\dot{b} = - \frac{2}{3}\, b \, \tr \DDkpt,\label{p15}
\end{align}
and
\begin{equation}
\xi=\left(\TT_\delta - [\TTel]_\delta \right) :\DD + \left(\mu (b - 3) - 2\sigma b\,\Delta b \right)\frac{\tr \DDkpt}{3} \quad \textrm{ with } \quad \xi \ge 0,
\label{p16}
\end{equation}
where
\begin{align*}
[\TTel]_\delta = - \sigma  (\nablax b \otimes \nablax b)_\delta.
\end{align*}
Requiring that
\begin{equation}
\begin{split}
\TT_\delta - [\TTel]_\delta & =  2\nu\, \DD \qquad \mbox{with $\nu>0$},\\
\mu (b - 3) - 2\sigma b\,\Delta b &= 2 \nu_1 \frac{\tr \DDkpt}{3} \qquad \mbox{with $\nu_1>0$},
\end{split}\label{p17}
\end{equation}
we obtain
\begin{equation}
\xi=2\nu |\DD|^2 + 2\nu_1 \frac{|\tr \DDkpt|^2}{9}.
\label{p16a}
\end{equation}
Referring then to \eqref{p15} we deduce that \eqref{p17} leads to
\begin{equation}
\begin{split}
\TT  = m\II + 2\nu \DD - \sigma ( \nabla b \otimes \nabla b)_\delta &= \phi \II + 2\nu \DD - \sigma ( \nabla b \otimes \nabla b), \\
\nu_1 \frac{\dot{b}}{b} + \mu (b - 3) - 2\sigma b\,\Delta b &= 0,
\end{split}\label{p18}
\end{equation}
and \eqref{p16a} implies that
\begin{equation}
\xi=2\nu |\DD|^2 + \frac{\nu_1}{2} \bigg|\frac{\dot{b}}{b}\bigg|^2.
\label{p16b}
\end{equation}
Note that $m$ and consequently $\phi:= m - \tfrac{\sigma}{3}|\nabla b|^2$ cannot be determined constitutively because of the incompressibility constraint. Consequently, it is an additional unknown quantity. In what follows we relabel $\phi$ by $-p$, as this symbol is used in most studies concerning incompressible fluids. To conclude, recalling also the definition of the material time derivative, the system of governing equations takes the form
\begin{align}
\divx \vv&=0, \label{p20}\\
\varrho\left(\partial_t \vv + \divx(\vv \otimes \vv)\right) - \divx \TT&=0, \label{p21} \\
\TT &=-p\,\II + 2\nu\DD - \sigma (\nablax b \otimes \nablax b), \label{p22} \\
\nu_1 \partial_t b + \nu_1 \divx (b\vv) + \mu (b^2 - 3b) - 2\sigma b^2 \Delta b  &= 0. \label{p23}
\end{align}
We wish to emphasize that despite the simplification $[\CCkpt]_{\delta} = \OO$ (see \eqref{ass1}), the reduced model \eqref{p20}--\eqref{p23} retains all of the material parameters of the original problem in the framework. It is also worth emphasizing that not only does the diffusion term appear in the equation for the mean normal elastic stress $b$, but also the Korteweg stress is present in the constitutive relation for the Cauchy stress.

Setting
\[ \tilde{b}:= \frac{b}{3} \quad \textrm{ and } \quad \tilde\sigma:= 9 \sigma \]
we can rewrite \eqref{p20}--\eqref{p23} as follows:
\begin{align}
\divx \vv&=0, \label{p20a}\\
\varrho(\partial_t \vv + \divx(\vv \otimes \vv)) - \divx \TT&=0, \label{p21a} \\
\TT &=-p\,\II + 2\nu\DD - \tilde \sigma (\nablax \tilde{b} \otimes \nablax \tilde{b}), \label{p22a} \\
\nu_1 \partial_t \tilde{b} + \nu_1 \divx (\tilde{b}\vv) + 3 \mu ({\tilde{b}}^2 - \tilde{b}) - 2\tilde\sigma {\tilde{b}}^2 \Delta \tilde{b}  &= 0, \label{p23a}
\end{align}
which, upon relabeling $b:=\tilde{b}$, $\sigma:=\tilde{\sigma}$ and $\mu:= 3\mu$ coincides with the system \eqref{eq10}--\eqref{eq12} that is analyzed in the next sections.

We conclude this section by showing that the choice of the constitutive relations for $\psi$ and $\xi$ directly yields the relevant \textit{a priori} estimates in a straightforward way. Note first that with the simplifications \eqref{ass1} and relabeling ($9\sigma$ by $\sigma$ and $3\mu$ by $\mu$) the constitutive relation \eqref{p10} reduces to
\begin{equation}
\varrho \psi = \varrho \psi_0(\varrho) + \frac{\mu}{2}(b-1 - \ln b) + \frac{\sigma}{2}|\nablax\, b|^2. \label{p10a}
\end{equation}
Next, by taking the scalar product of \eqref{p21a} with $\vv$ we obtain
\begin{equation} \label{apriori.1}
\partial_t \left(\varrho \frac{|\vv|^2}{2}\right) + \divx\bigg(\varrho \frac{|\vv|^2}{2} \, \vv\bigg) - \divx (\TT\vv) + \TT:\DD=0.
\end{equation}
Using the equation \eqref{p5} then leads to
\begin{equation} \label{apriori.2}
\partial_t \left(\varrho \frac{|\vv|^2}{2}\right) + \varrho\dot{\psi} + \xi + \divx\left(\varrho \frac{|\vv|^2}{2} \vv - \TT\vv + \qe - \theta \qn\right) = 0.
\end{equation}
Inserting \eqref{p16a} and \eqref{p10a} in \eqref{apriori.2}, and recalling that $\dot{\psi} = \partial_t\psi + (\vv\cdot \nabla) \psi$ and $\varrho\psi_0(\varrho) - \frac{\mu}{2}$ is constant, we arrive at
\begin{equation}\label{eq.geA}
\begin{split}
&\frac{1}{2}\partial_t \left(\varrho |\vv|^2 +\sigma  |\nabla b|^2+ \mu  (b-\ln b)\right)+ 2\nu |\DD|^2+\frac{\nu_1}{2} \bigg|\frac{\dot{b}}{b}\bigg|^2 \\
&\,\,+\divx\left( \frac{1}{2}\left(\varrho|\vv|^2+\mu (b- \ln b) +\sigma  |\nabla b|^2\right) \vv - \TT\, \vv + \sigma \dot{b}\nabla b  \right) =0.
\end{split}
\end{equation}
The next section starts with the derivation of  \eqref{eq.geA} from a PDE-analytic point of view.

\section{Function spaces and statement of the main result} \label{FS}

In this section, we provide the precise statement of the main theorem regarding the long-time existence of large-data weak solutions to the problem \eqref{eq10}--\eqref{data1}.

To motivate the choice of function spaces, we first briefly derive the basic energy estimate. For these formal calculations, we assume that all quantities are well defined; in particular, we assume that $b$ is strictly positive. Taking the scalar product of \eqref{eq10.5} with $\vv$ and using the incompressibility constraint \eqref{eq10}, we obtain the following identity (for the evolution of the kinetic energy):
\begin{equation*}
\frac12 \partial_t (\varrho|\vv|^2) + \frac12 \divx (\varrho \vv |\vv|^2 ) - \divx (\TT\, \vv) + \TT : \DD =0.
\end{equation*}
Employing further \eqref{eq11} and again \eqref{eq10} we obtain
\begin{equation}\label{eq.ke}
\frac12 \partial_t (\varrho |\vv|^2) + 2\nu |\DD|^2  =  \sigma (\nablax b \otimes \nablax b) : \DD- \frac12 \divx (\varrho \vv |\vv|^2 ) + \divx (\TT\, \vv) .
\end{equation}
Next, dividing \eqref{eq12} by $b^2$ and multiplying the result by $(\partial_t b + \divx (b\vv))$, we obtain
\begin{equation}\label{eq.bhelp}
\begin{aligned}
\frac{\nu_1}{b^2} \left|\partial_t b \right. & \left. + \divx (b\vv)\right|^2= (- \mu (1 - b^{-1}) + 2\sigma \Delta b)(\partial_t b + \divx (b\vv))\\
&=-\mu \left( \partial_t (b-\ln b) + \divx (\vv(b-\ln b))\right) -\sigma \partial_t |\nabla b|^2 + 2\sigma \divx ( \partial_t b \nabla b)\\
 &\quad+2\sigma \divx (\nabla b(\divx (b\vv)))-2\sigma \nabla b \cdot \nabla (\divx (b\vv)).
\end{aligned}
\end{equation}
Using the fact that $\divx \vv=0$ several times, the last term on the right-hand side can be rewritten as
$$
\begin{aligned}
\nabla b \cdot \nabla (\divx (b\vv))&= \nabla b \cdot \nabla (\vv \cdot \nabla b)\\ 
&=(\nabla b\otimes \nabla b) : \DD + \frac12 \divx (\vv |\nabla b|^2).
\end{aligned}
$$
Using this expression, \eqref{eq.bhelp} takes the form
\begin{equation}\label{eq.bhelp2}
\begin{aligned}
&\partial_t \left(\sigma  |\nabla b|^2+ \mu  (b-\ln b)\right)+\frac{\nu_1}{b^2} \left|\partial_t b + \divx (b\vv)\right|^2 =-2\sigma (\nabla b\otimes \nabla b) : \DD\\
&  + \divx\left( 2\sigma (\partial_tb+ \divx (b\vv))\nabla b - (\mu (b-\ln b) +\sigma  |\nabla b|^2) \vv\right).
\end{aligned}
\end{equation}
Finally, dividing  \eqref{eq.bhelp2} by 2 and adding the result to  \eqref{eq.ke},  we deduce the energy identity
\begin{equation}\label{eq.ge}
\begin{split}
&\frac{1}{2}\partial_t \left(\varrho|\vv|^2 +\sigma  |\nabla b|^2+ \mu  (b-\ln b)\right)+ 2\nu |\DD|^2+\frac{\nu_1}{2b^2} \left|\partial_t b + \divx (b\vv)\right|^2 \\
&\,\,=\divx\left(\TT\, \vv+ \sigma (\partial_tb+ \divx (b\vv))\nabla b - \frac{1}{2}\left(\varrho|\vv|^2+\mu (b-\ln b) +\sigma  |\nabla b|^2\right) \vv \right) .
\end{split}
\end{equation}
Hence, integration of the result over $\Omega\times (0,T)$, using integration by parts and the boundary conditions \eqref{eq.bc}, we deduce that
\begin{equation*}
\begin{split}
&\sup_{t\in (0,T)} \int_{\Omega} \left(\varrho|\vv|^2 +\sigma  |\nabla b|^2+ \mu  (b-\ln b)\right)\dd x + \int_Q 4\nu |\DD|^2+\frac{\nu_1}{b^2} \left|\partial_t b + \divx (b\vv)\right|^2\dd x \dd t \\
&\qquad \le \int_{\Omega} \left(\varrho|\vv_0|^2 +\sigma  |\nabla b_0|^2+ \mu  (b_0-\ln b_0)\right)\dd x.
\end{split}
\end{equation*}
This natural energy estimate thus motivates the minimal assumptions on the data $\vv_0$ and $b_0$ and simultaneously indicates the correct choices of function spaces for the problem. In addition, as it will be seen later, the unknown $b$ also satisfies a version of minimum/maximum principle and therefore we shall also require that $b_0$ is bounded and uniformly positive on $\Omega$.

Before we formulate the main theorem of the paper, we briefly recall the notations for the relevant function spaces. The standard Lebesgue and Sobolev spaces are denoted by $L^p(\Omega)$ and $W^{k,p}(\Omega)$, respectively, and are equipped with the usual norms denoted by $\|\cdot\|_p$ and $\|\cdot\|_{k,p}$, respectively. Since we also need to deal with vector- and tensor-valued functions, we will emphasize this character by writing  $L^p(\Omega)^d$ or $L^p(\Omega)^{d\times d}$ if necessary. Moreover, since the velocity is a solenoidal function, we define
$$
\begin{aligned}
{\color{black} L^2_{0,\divx}}&:= \overline{\left\{\vv \in \mathcal{C}^{\infty}_0(\Omega; \mathbb{R}^d): \; \divx \vv =0\right\}}^{\|\cdot \|_2},\\
{\color{black} W^{1,2}_{0,\divx}}&:= \overline{\left\{\vv \in \mathcal{C}^{\infty}_0(\Omega; \mathbb{R}^d): \; \divx \vv =0\right\}}^{\|\cdot \|_{1,2}}.
\end{aligned}
$$

Having defined the proper function spaces, we notice that the assumptions on the data, see \eqref{data1}, are a natural consequence of \eqref{eq.ge}.
We can now formulate our main theorem.
\begin{theorem}\label{T:main}
Let $\Omega \subset \mathbb{R}^d$ be a Lipschitz domain and $T>0$. Assume that $\varrho$, $\nu$, $\nu_1$, $\mu$ and $\sigma$ are positive constants and that $(\vv_0, b_0)$ satisfies \eqref{data1}. Then, there exists a couple $(\vv,b)$ such that
\begin{align}
\vv &\in L^{\infty}(0,T; L^2_{0,\divx}) \cap L^2(0,T; W^{1,2}_{0,\divx}),\label{S1}\\
b&\in L^{\infty}(0,T; W^{1,2}(\Omega)), \quad b, b^{-1} \in L^{\infty}(Q), \label{S2}\\
\partial_t \vv&\in L^1(0,T; (W^{1,2}_{0,\divx}\cap W^{d+1,2}(\Omega)^d)^*)\label{S3},\\
\partial_t b & \in L^1(Q)\label{S4}
\end{align}
and
\begin{equation}\label{S5}
\begin{split}
(\partial_t b + \divx(b\vv)) &\in L^2(Q),\\
\Delta b &\in L^2(Q),
\end{split}
\end{equation}
and this couple solves \eqref{eq10}--\eqref{eq.ic} in the following sense: for almost all $t\in (0,T)$ and all $\vw \in W^{1,2}_{0,\divx}\cap W^{d+1,2}(\Omega)^d$ we have (recall that $\DD$ denotes the symmetric part of $\nabla \vv$)
\begin{equation}
\langle \partial_t (\varrho\vv), \vw \rangle + \int_{\Omega} \left(2\nu \DD - \sigma \nabla b \otimes \nabla b - \varrho \vv \otimes \vv\right) : \nabla \vw \dd x =0; \label{Tw1}
\end{equation}
for almost all $t\in (0,T)$ and all $w\in W^{1,2}(\Omega)$ we have
\begin{equation}
\int_{\Omega} \left(\frac{\nu_1 (\partial_t b +  \divx (b\vv))}{b^2} + \mu \left(1 - \frac{1}{b}\right)\right)w +  2\sigma \nabla  b \cdot \nabla w \dd x = 0, \label{Tw2}
\end{equation}
and the initial data are attained in the following sense:
\begin{equation}
\lim_{t\to 0_+} (\|\vv(t)-\vv_0\|_2+\|b(t)-b_0\|_{1,2}) =0.\label{T:ic}
\end{equation}
Moreover, the following energy inequality holds for all $t\in (0,T)$:
\begin{equation}
\begin{split}
&\frac12 \int_{\Omega} \left(\varrho|\vv(t)|^2 +\sigma  |\nabla b(t)|^2+ \mu  (b(t)-\ln b(t))\right)\dd x \\
&\quad + \int_0^t\int_{\Omega} 2\nu |\DD|^2+\frac{\nu_1}{2b^2} \left|\partial_t b + \divx (b\vv)\right|^2\dd x \dd \tau \\
&\qquad \le \frac12\int_{\Omega} \left(\varrho|\vv_0|^2 +\sigma  |\nabla b_0|^2+ \mu  (b_0-\ln b_0)\right)\dd x.
\end{split}\label{T:EI}
\end{equation}
\end{theorem}
Notice here that \eqref{eq10} is automatically satisfied since {\color{black} $\vv\in L^\infty(0,T;L^2_{0,\divx})$}; the identity \eqref{Tw1} represents the weak formulation of \eqref{eq10.5}, \eqref{eq11} and \eqref{eq.bc}$_1$, where we have omitted the presence of the pressure $p$ by choosing divergence-free test functions. Finally, \eqref{Tw2} is a weak formulation of \eqref{eq12} and \eqref{eq.bc}$_2$, which is obtained by division of \eqref{eq12} by $b^2$, multiplying the result by a test function $w$, integrating over $\Omega$ and in the elliptic term using integration by parts together with   \eqref{eq.bc}$_2$. Observe also that because of the assumed regularity \eqref{S2} and \eqref{S5}, we could simply say that \eqref{eq12} holds almost everywhere in $Q$, but we would lose the information \eqref{eq.bc}$_2$, which is encoded in \eqref{Tw2}. Nevertheless, \eqref{eq12} can be directly obtained from \eqref{Tw2} by setting ${\color{black} w:=b^2 \varphi}$ with an arbitrary $\varphi \in \mathcal{C}^{\infty}_0(\Omega)$ and noting that (here we use the fact that $\Delta b \in L^2(\Omega)$ and integration by parts)
$$
\begin{aligned}
\int_{\Omega} \nabla b \cdot \nabla (b^2 \varphi)\dd x&= \int_{\Omega} 2b|\nabla b|^2 \varphi + b^2 \nabla b \cdot \nabla \varphi\dd x=\int_{\Omega} 2b|\nabla b|^2 \varphi - \divx ( b^2 \nabla b) \varphi\dd x\\
&=-\int_{\Omega} b^2\Delta b \varphi\dd x.
\end{aligned}
$$

\section{Proof of Theorem~\ref{T:main}} \label{Proof}
For the sake of simplicity we set the parameters $\nu_1$ and $\mu$ equal to 2 and $\varrho$, $\nu$ and $\sigma$ equal to 1, and we also divide \eqref{eq12} by 2 in order to avoid the presence of all constants. Note that such a simplification does not change the assertion of the theorem, as the proof can be straightforwardly extended to any other values of these parameters thanks to their positivity. We will also frequently use the symbol $C$ to denote a generic positive constant, whose value may change from line to line, and can depend only on the data. In instances when $C$ also depends on the indices of the
approximating problem appearing in the proof, this will be clearly indicated.

\subsection{Galerkin approximation}

First, we introduce the following cut-off function
$$
T_{n}(s):=\min \left\{n, \max\{n^{-1},s \} \right\} \textrm{ for } s\in \mathbb{R}.
$$
Next, we find a basis $\{\vw_n\}_{n=1}^{\infty}$ of $W^{1,2}_{0,\divx} \cap W^{d+1,2}(\Omega;\mathbb{R}^d)$, which is orthogonal in $L^2_{0, \divx}$ and a basis {\color{black} $\{w_n\}_{n=1}^{\infty}$ of $W^{1,2}(\Omega)$, which is orthonormal in $L^2( \Omega)$ and orthogonal in $W^{1,2}(\Omega)$. The existence of the latter follows from a version of the Hilbert--Schmidt theorem (cf. Lemma 5.1 in \cite{FS2012}), according to which $w_n$, $n\in \mathbb{N}$, can be taken to be the eigenfunctions of the operator $w \mapsto -\Delta w + w$ subject to a homogeneous Neumann boundary condition
on $\partial \Omega$. By Theorem 1.2 in \cite{Geng} then $w_n \in W^{1,p}(\Omega)$ for all $p \in (\frac{2d}{d+1} - \epsilon, \frac{2d}{d-1} + \epsilon)$, $d \in \{2,3\}$. Hence, by Morrey's embedding theorem $w_n \in C^{0,\alpha}(\overline{\Omega})$, whereby $\Delta w_n \in C^{0,\alpha}(\overline{\Omega})$, for some $\alpha \in (0,\frac{1}{2})$.} Let, for arbitrary $n \in \mathbb{N}$ and $\ell\in \mathbb{N}$, the symbols $\vec{P}^n$ and $P^{\ell}$ denote the orthogonal (in $L^2$) projections onto $\textrm{span}\{\vw_i\}_{i=1}^{n}$ and $\textrm{span}\{w_i\}_{i=1}^{\ell}$, respectively. {\color{black}
Observe also that $\|P^\ell(w)\|_{1,2} \leq \|w\|_{1,2}$ for all $w \in W^{1,2}(\Omega)$ and all $\ell \in \mathbb{N}$}. Finally, we consider
$$
\vv^{n,\ell}(t,x):= \sum_{i=1}^n \alpha_{i}^{n,\ell}(t) \vw_i(x),\qquad b^{n,\ell}(t,x):= \sum_{i=1}^{\ell} \beta_i^{n,\ell}(t) w_i(x),
$$
satisfying
$$
\vv^{n,\ell}(0,x)=\vv_0^n(x):= \vec{P}^n(\vv_0)(x), \qquad b^{n,\ell}(0,x)=b_0^{\ell}(x):= P^{\ell}(b_0)(x),
$$
and solving, over  $(0,T)$,  the following system\footnote{Notice that \eqref{eq.G2} is an approximation of \eqref{eq12} after division by $b^2$, i.e., an approximation of \eqref{Tw2}.} of ordinary differential equations:
\begin{align}
\begin{aligned}
\int_{\Omega}\partial_t \vv^{n,\ell} \cdot \vw_i - \vv^{n,\ell} \otimes \vv^{n,\ell} : \nabla \vw_i + \bS^{n,\ell} : \nabla \vw_i &\dd x = 0,
\\&\qquad i=1,\ldots, n, \end{aligned} \label{eq.G1} \\
\begin{aligned}
\int_{\Omega} \frac{\partial_t b^{n,\ell} w_j}{(T_n(b^{n,\ell}))^2} &+ \frac{\nabla b^{n,\ell}\cdot \vv^{n,\ell}}{(T_n(b^{n,\ell}))^2}w_j
+ (1 - (T_n(b^{n,\ell}))^{-1})w_j \dd x \\ &+ \int_{\Omega}  \nabla b^{n,\ell}\cdot \nabla w_j \dd x  = 0, \qquad j=1,\ldots, \ell, \end{aligned} \label{eq.G2}
\end{align}
where $\bS^{n,\ell}$ is defined a.e. in $Q$ by
\begin{equation}
\bS^{n,\ell}:=  2\DD^{n,\ell} - (\nablax b^{n,\ell} \otimes \nablax b^{n,\ell}), \quad \DD^{n,\ell}:=\frac{1}{2}\big(\nabla \vv^{n,\ell}+(\nabla \vv^{n,\ell})^{\textrm{T}}\big). \label{eq.G3}
\end{equation}
To prove the existence of a solution to the above problem is not a difficult task but since it requires more than invoking standard Carath\'{e}odory theory, we provide a sketch of the proof in the Appendix for the sake of completeness.

\subsection{The limit $\ell \to \infty$}

We start this part by establishing the estimates that will be independent of $\ell$ but can still be dependent on $n$. Nevertheless, these estimates will suffice for the purpose of letting $\ell \to \infty$.

First, we define the following primitive function:
$$
\Theta_n(s):= \int_0^s \frac{s}{(T_n(s))^2} \dd s.
$$
It follows directly from the properties of the function $T_n$ that
\begin{equation}
n^{-2}s^2\le 2\Theta_n(s)\le n^2s^2. \label{Tkb}
\end{equation}
Then we multiply the $j$-th equation in \eqref{eq.G2} by $\beta^{n,\ell}_j(t)$ and sum the resulting identities over $j=1,\ldots, \ell$ to obtain the following equality:
\begin{equation}
\label{eq.ID1}
\int_{\Omega} \frac{\partial_t b^{n,\ell} b^{n,\ell}}{(T_n(b^{n,\ell}))^2} + \frac{b^{n,\ell}\nabla b^{n,\ell}\cdot \vv^{n,\ell}}{(T_n(b^{n,\ell}))^2} + (1 - (T_n(b^{n,\ell}))^{-1})b^{n,\ell} + |\nabla b^{n,\ell}|^2 \dd x  = 0.
\end{equation}
Using the definition of $\Theta_n$ and the facts that $\divx \vv^{n,\ell}=0$ and $\vv^{n,\ell}$ vanishes on $\partial \Omega$, the second term vanishes and the above identity reduces to the following inequality (using also the properties of $T_n$ and $\Theta_n$):
$$
\frac{\dd}{\dd t}\int_{\Omega} \Theta_n(b^{n,\ell}(t))\dd x  + \|\nabla b^{n,\ell}(t)\|_2^2 \le C(n)\left(1+\int_{\Omega}\Theta_n(b^{n,\ell})\dd x \right).
$$
Thus, using Gronwall's lemma, we have the following $\ell$-independent bound:
\begin{equation}
\sup_{t\in (0,T)}\|b^{n,\ell}(t)\|_2^2 + \int_0^T \|b^{n,\ell}(t)\|_{1,2}^2 \dd t \le C(n) (1+\|b_0\|_2^2).\label{eq.est1}
\end{equation}

Next, we focus on $\ell$-independent bounds on $\vv^{n,\ell}$. Since $n$ is fixed, we first recall that (in an $n$-dimensional space)
\begin{equation}
\|\vv^{n,\ell}(t)\|_{1,\infty}\le C(n) \|\vv^{n,\ell}(t)\|_2.\label{eq.eq}
\end{equation}
Furthermore, multiplying the $i$-th equation in \eqref{eq.G1} by $\alpha_i^{n,\ell}$ and summing the result over $i=1,\ldots, n$, using the definition of $\bS^{n,\ell}$, we get that (using again the fact that $\divx \vv^{n,\ell}=0$ and so the convective term vanishes)
\begin{equation}\label{eq.ID3}
\begin{split}
\frac12 \frac{\dd}{\dd t} \|\vv^{n,\ell}(t)\|_2^2 &+ 2\|\DD^{n,\ell}(t)\|_2^2 = \int_{\Omega}\nabla b^{n,\ell}(t) \otimes \nabla b^{n,\ell}(t) : \DD^{n,\ell}(t)\dd x \\
&\le \|\nabla b^{n,\ell}(t)\|_2^2 \|\vv^{n,\ell}(t)\|_{1,\infty} \\
&\le C(n)\|\nabla b^{n,\ell}(t)\|_2^2 \|\vv^{n,\ell}(t)\|_{2} \\
&\le C(n)\|\nabla b^{n,\ell}(t)\|_2^2 + C(n) \|\nabla b^{n,\ell}(t)\|_2^2 \|\vv^{n,\ell}(t)\|^{2}_{2},
\end{split}
\end{equation}
where we have also used \eqref{eq.eq}. Therefore, thanks to the estimate \eqref{eq.est1}, Gronwall's lemma implies that
\begin{equation}
\begin{split}
\sup_{t\in (0,T)}\|\vv^{n,\ell}(t)\|_{1,\infty}^2&\le C(n)\sup_{t\in (0,T)}|\alpha^{n,\ell}(t)|^2 \\
&\le C(n) \sup_{t\in (0,T)}\|\vv^{n,\ell}(t)\|_2^2 \le C(n,\|\vv_0\|_2, \|b_0\|_2).\label{eq.est2}
\end{split}
\end{equation}
Finally, we derive uniform bounds on time derivatives. To do so, we multiply the $j$-th equation in \eqref{eq.G2} by $\partial_t \beta_j^{n,\ell}$ to deduce the following identity:
\begin{equation}\label{eq.ID2}
\begin{aligned}
&\frac{1}{2}\frac{\dd}{\dd t} \int_{\Omega} |\nabla b^{n,\ell}|^2\dd x+ \int_{\Omega} \frac{|\partial_t b^{n,\ell}|^2}{(T_n(b^{n,\ell}))^2}\dd x \\
&\qquad = \int_{\Omega} ((T_n(b^{n,\ell}))^{-1}-1)\,\partial_t b^{n,\ell} -\frac{\partial_t b^{n,\ell}\nabla b^{n,\ell}\cdot \vv^{n,\ell}}{(T_n(b^{n,\ell}))^2}\dd x.
\end{aligned}
\end{equation}
Hence, by employing the estimate \eqref{eq.est2}, Young's inequality and the properties of $T_n$, we get
$$
\frac{\dd}{\dd t} \|\nabla b^{n,\ell}(t)\|_2^2 + \|\partial_t b^{n,\ell}(t)\|_2^2 \le C(n) (1+ \|\nabla b^{n,\ell}(t)\|_2^2),
$$
so that by applying Gronwall's lemma again we have
\begin{equation}\label{eq.est3}
\sup_{t\in (0,T)} \|b^{n,\ell}(t)\|_{1,2}^2 + \int_0^T \|\partial_t b^{n,\ell}(t)\|_2^2 \dd t \le C(n,\|b_0\|_{1,2}, \|\vv_0\|_2).
\end{equation}
Finally, using the last inequality and also \eqref{eq.est2}, we see from \eqref{eq.G1} that
\begin{equation}\label{eq.est4}
\sup_{t\in (0,T)} \left|\frac{\dd}{\dd t} \alpha^{n,\ell}(t)\right|\le C(n).
\end{equation}

Having established bounds that are uniform with respect to $\ell$, we can now focus on taking the limit $\ell \to \infty$ in \eqref{eq.G1}, \eqref{eq.G2}. Indeed, using \eqref{eq.est2} and \eqref{eq.est4}, there exists a subsequence that we do not relabel such that
\begin{align}
\alpha^{n,\ell}&\rightharpoonup^* \alpha^n &&\textrm{weakly$^*$ in } W^{1,\infty}(0,T; \mathbb{R}^n),\label{c1}\\
\alpha^{n,\ell}&\to \alpha^n &&\textrm{strongly in }\mathcal{C}([0,T]; \mathbb{R}^n).\label{c2}
\end{align}
Thus, it follows from the definition of $\vv^{n,\ell}$ that
\begin{align}
\vv^{n,\ell}&\rightharpoonup^* \vv^n &&\textrm{weakly$^*$ in } W^{1,\infty}(0,T; W^{1,\infty}(\Omega)^d),\label{c1b}\\
\vv^{n,\ell}&\to \vv^n &&\textrm{strongly in }\mathcal{C}([0,T]; W^{1,\infty}(\Omega)^d),\label{c2b}
\end{align}
where
$$
\vv^n(t,x)=\sum_{i=1}^n \alpha_i^n(t)\vw_i(x) \qquad \textrm{and} \qquad \vv^n(0,x)=\vec{P}^n(\vv_0)(x).
$$
Similarly, using \eqref{eq.est3} and compact embedding, we find a subsequence (that is again labeled in the same way) such that
\begin{align}
b^{n,\ell}&\rightharpoonup^* b^n &&\textrm{weakly$^*$ in } L^{\infty}(0,T; W^{1,2}(\Omega)),\label{c3}\\
b^{n,\ell}&\rightharpoonup b^n &&\textrm{weakly in }W^{1,2}(0,T; L^2(\Omega)),\label{c4}\\
b^{n,\ell}&\to b^n &&\textrm{strongly in }L^2(0,T; L^2(\Omega)) \textrm{ and a.e. in } Q.\label{c5}
\end{align}
In addition, it follows directly from the construction that $b^n(0)=b_0$.

The convergence results \eqref{c2b} and \eqref{c3}--\eqref{c5} allow us to let $\ell \to \infty$ in \eqref{eq.G2} directly, so we conclude that for almost all $t\in (0,T)$ and all $w\in W^{1,2}(\Omega)$ we have
\begin{equation}
\int_{\Omega} \frac{\partial_t b^{n} w}{(T_n(b^{n}))^2} + \frac{\nabla b^{n}\cdot \vv^{n}}{(T_n(b^{n}))^2}w + (1 - (T_n(b^{n}))^{-1})w + \nabla b^{n}\cdot \nabla w \dd x  = 0.\label{eq.Gn2}
\end{equation}
In particular, setting $w:=b^n(t)$ in the above identity and integrating the result over $(0,T)$ we obtain
\begin{equation}
\int_Q|\nabla b^{n}|^2\dd x \dd t=-\int_{Q} \frac{\partial_t b^{n} b^n}{(T_n(b^{n}))^2} + \frac{\nabla b^{n}\cdot \vv^{n}}{(T_n(b^{n}))^2}b^n + (1 - (T_n(b^{n}))^{-1})b^n\dd x \dd t.\label{eq.Gn2a}
\end{equation}
Furthermore, integration of \eqref{eq.ID1} over $(0,T)$ and the use of \eqref{c2b} and \eqref{c3}--\eqref{c5} lead to
\begin{equation}
\label{eq.ID1n}
\begin{split}
&\lim_{\ell \to \infty} \int_Q |\nabla b^{n,\ell}|^2 \dd x \dd t\\
&=-\lim_{\ell \to \infty} \int_{Q} \frac{\partial_t b^{n,\ell} b^{n,\ell}}{(T_n(b^{n,\ell}))^2} + \frac{b^{n,\ell}\nabla b^{n,\ell}\cdot \vv^{n,\ell}}{(T_n(b^{n,\ell}))^2} + (1 - (T_n(b^{n,\ell}))^{-1})b^{n,\ell}\dd x \dd t\\
&=-\int_{Q} \frac{\partial_t b^{n} b^n}{(T_n(b^{n}))^2} + \frac{\nabla b^{n}\cdot \vv^{n}}{(T_n(b^{n}))^2}b^n + (1 - (T_n(b^{n}))^{-1})b^n\dd x \dd t\\
&=\int_Q |\nabla b^{n}|^2 \dd x \dd t,
\end{split}
\end{equation}
where we have used \eqref{eq.Gn2a} for the last equality. Consequently, combining \eqref{c3} and \eqref{eq.ID1n} we deduce
that
\begin{align}
b^{n,\ell}&\to b^n &&\textrm{strongly in } L^{2}(0,T; W^{1,2}(\Omega)).\label{c6}
\end{align}
Finally, this strong convergence result together with \eqref{c1} and \eqref{c2} allows us to let $\ell \to \infty$ in \eqref{eq.G1} and \eqref{eq.G3} to get, {\color{black} for a.e. $t \in (0,T)$,}
\begin{align}
\int_{\Omega}\partial_t \vv^{n} \cdot \vw_i - \vv^{n} \otimes \vv^{n} \cdot \nabla \vw_i + \bS^{n}\cdot \nabla \vw_i\dd x&=0, \qquad i=1,\ldots, n,\label{eq.G1n}
\end{align}
where {\color{black} $\{\vw_i\}_{i=1}^n$ is the set of Galerkin basis functions defined at the start of the section; furthermore, } $\bS^{n}$ is defined by
\begin{equation}
\bS^{n}=  2\DD^{n} - (\nablax b^{n} \otimes \nablax b^{n}) \quad \mbox{and} \quad \DD^{n}=\frac{1}{2}\big(\nabla \vv^{n}+(\nabla \vv^{n})^{\textrm{T}}\big). \label{eq.G3n}
\end{equation}
{\color{black} We note here that the precise passage to the limit in \eqref{eq.G1} as $\ell \rightarrow \infty$, leading to \eqref{eq.G1n}, requires multiplication of \eqref{eq.G1} by any $\varphi \in C^\infty_0(0,T)$ and integration of the resulting equality over $t \in (0,T)$, followed by passage to the limit therein as $\ell \rightarrow \infty$, which then results in \eqref{eq.G1n} multiplied by $\varphi \in C^\infty_0(0,T)$ and integrated over $t \in (0,T)$. Thanks to the fact that each term under the integral sign in the limiting equality is locally
integrable with respect to $t \in (0,T)$ and the equality holds for all $\varphi \in C^\infty_0(0,T)$, the fundamental theorem of the calculus of variations (Du Bois Reymond's lemma) yields the equality \eqref{eq.G1n} for a.e. $t \in (0,T)$, as stated.}

\subsection{Estimates independent of $n$}
We first focus on minimum and maximum principles for $b^n$. To this end, we define
$$
b_{\max}:=\max \left\{1,\|b_0\|_{\infty}\right\}, \qquad \frac{1}{b_{\min}}:=\max\left\{1,\left\|\frac{1}{b_0}\right\|_{\infty}\right\},
$$
and we also recall the notations $s_+:=\max\{0,s\}$ and $s_{-}:=\min\{0,s\}$. From now on, we shall assume that $n$ is sufficiently large; more precisely, we shall suppose that $n>b_{\max}$ and $n>(b_{\min})^{-1}$.

Then, we set $w:=(b^n-b_{\max})_+$ in \eqref{eq.Gn2} and using the facts that $b_{\max}\ge 1$ and $n>1$, we have
\begin{equation}
\begin{split}
&\int_{\Omega} \frac{\partial_t b^{n} (b^n-b_{\max})_+}{(T_n(b^{n}))^2} + \frac{\nabla b^{n}\cdot \vv^{n}}{(T_n(b^{n}))^2}(b^n-b_{\max})_+ \dd x \\
&=-\int_{\Omega} (1 - (T_n(b^{n}))^{-1})(b^n-b_{\max})_+ + \nabla b^{n}\cdot \nabla (b^n-b_{\max})_+ \dd x\\
&=-\int_{b^n>b_{\max}} \frac{T_n(b_n)-1}{T_n(b^{n})}(b^n-b_{\max})_+ + |\nabla b^{n}|^2\dd x\le 0.\label{eq.Gnmax}
\end{split}
\end{equation}
Hence, defining
$$
\Gamma_+(t):=\int_{b_{\max}}^t \frac{(s-b_{\max})_+}{(T_n(s))^2}\dd s \ge 0,
$$
we obtain with the help of \eqref{eq.Gnmax} that
\begin{equation}
\begin{split}
0&\ge \int_{\Omega} \frac{\partial_t b^{n} (b^n-b_{\max})_+}{(T_n(b^{n}))^2} + \frac{\nabla b^{n}\cdot \vv^{n}}{(T_n(b^{n}))^2}(b^n-b_{\max})_+ \dd x \\
&=\int_{\Omega} \partial_t \Gamma_{+} (b^n) + \nabla \Gamma_{+}(b^n) \cdot \vv^n \dd x=\frac{\dd}{\dd t}\int_{\Omega} \Gamma_{+} (b^n)\dd x,\label{eq.Gnmax2}
\end{split}
\end{equation}
where for the last equality we have used integration by parts in the second term and the fact that $\vv^n\in W^{1,2}_{0,\divx}$. Since $b_0\le b_{\max}$ a.e. in $\Omega$, the last inequality implies (using the positivity of $\Gamma_+(s)$ whenever $s>b_{\max}$) that
\begin{equation}
b^n\le b_{\max} \textrm{ a.e. in }Q.\label{maxp}
\end{equation}

Similarly, we deduce the minimum principle for $b^n$.  We set $w:=(b^n-b_{\min})_-$ in \eqref{eq.Gn2} and see that (note that $b_{\min}\le 1$ and $n>1$)
\begin{equation}
\begin{split}
&\int_{\Omega} \frac{\partial_t b^{n} (b^n-b_{\min})_-}{(T_n(b^{n}))^2} + \frac{\nabla b^{n}\cdot \vv^{n}}{(T_n(b^{n}))^2}(b^n-b_{\min})_- \dd x\\
&=-\int_{\Omega} (1 - (T_n(b^{n}))^{-1})(b^n-b_{\min})_- + \nabla b^{n}\cdot \nabla (b^n-b_{\min})_- \dd x\\
&=-\int_{b^n<b_{\min}} \frac{T_n(b_n)-1}{T_n(b^{n})}(b^n-b_{\min})_- + |\nabla b^{n}|^2\dd x\le 0.\label{eq.Gnmin}
\end{split}
\end{equation}
Hence, defining
$$
\Gamma_-(t):=-\int^{b_{\min}}_t \frac{(s-b_{\min})_-}{(T_n(s))^2}\dd s \ge 0,
$$
we obtain with the help of \eqref{eq.Gnmin} that
\begin{equation}
\begin{split}
0&\ge \int_{\Omega} \partial_t \Gamma_{-} (b^n) + \nabla \Gamma_{-}(b^n) \cdot \vv^n \dd x=\frac{\dd}{\dd t}\int_{\Omega} \Gamma_{-} (b^n)\dd x.\label{eq.Gnmin2}
\end{split}
\end{equation}
Thus, the last inequality implies (using the positivity of $\Gamma_-(s)$ whenever $s<b_{\min}$) that
\begin{equation}
b^n\ge b_{\min} \textrm{ a.e. in }Q.\label{minp}
\end{equation}

We continue with further uniform estimates. Having \eqref{maxp} and \eqref{minp}, we see that for sufficiently large $n$ there holds $T_n(b^n)=b^n$, and \eqref{eq.Gn2} then reduces to
\begin{equation}
\int_{\Omega} \frac{\partial_t b^{n} w}{(b^{n})^2} + \frac{\nabla b^{n}\cdot \vv^{n}}{(b^{n})^2}w + (1 - (b^{n})^{-1})w + \nabla b^{n}\cdot \nabla w \dd x  = 0\label{eq.Gn2m}
\end{equation}
for all $w\in W^{1,2}(\Omega)$. However, using \eqref{c3} and \eqref{c4}, we see from the above identity that, for each (large) $n$,  $\Delta b^n \in L^2(0,T; L^2(\Omega))$, and consequently, using interior elliptic regularity we have
\begin{equation}
b^n \in L^2(0,T; W^{2,2}(\tilde{\Omega})) \qquad \textrm{ for all } \tilde{\Omega} \subset \subset \Omega. \label{elip}
\end{equation}
In addition, using \eqref{eq.Gn2m}, almost everywhere in $Q$ we have
\begin{equation}
\frac{\partial_t b^{n} }{(b^{n})^2} + \frac{\nabla b^{n}\cdot \vv^{n}}{(b^{n})^2} + (1 - (b^{n})^{-1}) -\Delta b^{n}  = 0.\label{eq.Gn2mae}
\end{equation}

Next, we would like to mimic the procedure that led to the identity \eqref{eq.ge}. This in particular means that we would like to multiply \eqref{eq.Gn2mae} by $\partial_t b^n + \nabla b^n \cdot \vv^n$, integrate over $\Omega$, and then use integration by parts in the elliptic term. However, for such a procedure, we do not have sufficient regularity of $b^n$ and therefore we need to proceed more carefully. Thus, instead of testing \eqref{eq.Gn2mae} by $\partial_t b^n$, we pass to the limit in \eqref{eq.ID2}. Indeed, integration of \eqref{eq.ID2} over $(0,t)$ leads to
\begin{equation*}
\begin{aligned}
& \frac{1}{2}\|\nabla b^{n,\ell}(t)\|^2_2+ \int_0^t\int_{\Omega} \frac{|\partial_t b^{n,\ell}|^2}{(T_n(b^{n,\ell}))^2}\dd x \dd s\\
&= \int_0^t \int_{\Omega} ((T_n(b^{n,\ell}))^{-1}-1)\partial_t b^{n,\ell} -\frac{\partial_t b^{n,\ell}\nabla b^{n,\ell}\cdot \vv^{n,\ell}}{(T_n(b^{n,\ell}))^2}\dd x\dd s+{\color{black} \frac{1}{2} \|\nabla P^\ell(b_0)\|_2^2}.
\end{aligned}
\end{equation*}
Then, with the help of \eqref{c2b}, \eqref{c3}--\eqref{c5}, \eqref{c6}, the fact that $T_n(b^n)=b^n$ and the weak lower semicontinuity of norms, we obtain that for all\footnote{In fact we first obtain \eqref{eq.ID2int} only for almost all $t\in (0,T)$. However, thanks to the continuity of $b^n$ in the weak topology, we can extend the result onto the whole time interval $(0,T)$.} $t\in (0,T)$ one has the following inequality:
\begin{equation}\label{eq.ID2int}
\begin{aligned}
& \frac{1}{2}\|\nabla b^{n}(t)\|^2_2 + \int_0^t\int_{\Omega} \frac{|\partial_t b^{n}|^2}{(b^{n})^2}\dd x \dd s\\
&\quad \le \int_0^t \int_{\Omega} ((b^{n})^{-1}-1)\partial_t b^{n} -\frac{\partial_t b^{n}\nabla b^{n}\cdot \vv^{n}}{(b^{n})^2}\dd x\dd s+ {\color{black} \frac{1}{2}  \|b_0\|_{1,2}^2}\\
&\quad = -\int_0^t \int_{\Omega}\frac{\partial_t b^{n}\nabla b^{n}\cdot \vv^{n}}{(b^{n})^2}\dd x\dd s+\int_{\Omega} \ln b^n(t)-\ln b_0 -b^n(t)+b_0\dd x+ {\color{black}\frac{1}{2}\|b_0\|_{1,2}^2}.
\end{aligned}
\end{equation}

{\color{black}
We note in particular that the second term on the left-hand side of \eqref{eq.ID2int} arises by letting $\ell \rightarrow \infty$ in the second term on the left-hand side of the displayed equality preceding \eqref{eq.ID2int}, by applying the weak lower semicontinuity result stated in Theorem 3.23 on p.96 of \cite{dacorogna}, thanks to the strong convergence $b^{n,\ell} \to b^n$ in $L^2(0,T; L^2(\Omega))$ as $\ell \rightarrow \infty$ (cf. \eqref{c5}),
the weak convergence result $\partial_t b^{n,\ell} \rightharpoonup \partial_t b^n$ in  $L^2(0,T; L^2(\Omega))$ as $\ell \rightarrow \infty$ (cf. \eqref{c4}), and the convexity of the function $\xi \in \mathbb{R} \mapsto |\xi|^2$ featuring in the numerator of the integrand in the term concerned.}

Finally, we multiply \eqref{eq.Gn2mae} by $\nabla b^n \cdot \vv^n$ and integrate the result over $\Omega$ to get (notice here that thanks to the \textit{a~priori} estimates, in particular thanks to the fact that $\Delta b^n \in L^2(\Omega)$ for almost all $t\in (0,T)$, such a procedure is rigorous):
\begin{equation}
\int_{\Omega}\frac{|\nabla b^{n}\cdot \vv^{n}|^2}{(b^{n})^2} + \frac{\partial_t b^{n} \nabla b^{n}\cdot \vv^{n} }{(b^{n})^2}\dd x = \int_{\Omega} ((b^{n})^{-1}-1)\nabla b^{n}\cdot \vv^{n} +\Delta b^{n}\nabla b^{n}\cdot \vv^{n}\dd x.\label{eq.Gn2maes}
\end{equation}
First, we see that since $\divx \vv^n=0$ and $\vv^n$ has zero trace the first integral on the right-hand side can be evaluated with the help of integration by parts as follows:
\begin{equation}
\int_{\Omega} ((b^{n})^{-1}-1)\nabla b^{n}\cdot \vv^{n}\dd x=\int_{\Omega} \nabla (\ln b^n - b^n)\cdot \vv^{n}\dd x =0. \label{zero1}
\end{equation}
Next, we focus on the last term in \eqref{eq.Gn2maes}. We would like to integrate by parts, but since we do not know that $b^n$ has second derivatives integrable up to the boundary, we have to proceed more carefully. However, since the integral is well defined, we know that
$$
\int_{\Omega} \Delta b^{n}\,\nabla b^{n}\cdot \vv^{n}\dd x = \lim_{\varepsilon \to 0_+} \int_{\Omega}\Delta b^{n}\,\nabla b^{n}\cdot \vv^{n}\tau_{\varepsilon}\dd x,
$$
where $\tau_{\varepsilon}$ is smooth and satisfies $|\tau_{\varepsilon}|\le 1$ and
$$
\tau_{\varepsilon}(x)=\left\{
\begin{aligned}
&0 &&\textrm{ if } \textrm{dist }(x,\partial \Omega)\le \varepsilon,\\
&1 &&\textrm{ if } \textrm{dist }(x,\partial \Omega)\ge 2\varepsilon.
\end{aligned}
\right.
$$
Furthermore, since $\Omega$ is Lipschitz and $\vv^n \in W^{1,\infty}_{0,\divx}$, we can construct such a $\tau_{\varepsilon}$ so that also
\begin{equation}\label{tau}
|\vv^n||\nabla \tau_{\varepsilon}| \le C (\vv^n,\Omega),
\end{equation}
with $C$ independent of $\varepsilon$. Then, using integration by parts {\color{black} recalling again that $\Omega$ is Lipschitz, using} the fact that $\tau_{\varepsilon}$ is compactly supported, the fact that $\divx \vv^n= 0$, the estimate \eqref{tau} and the regularity \eqref{elip}, we get
$$
\begin{aligned}
&\int_{\Omega}\Delta b^{n}\nabla b^{n}\cdot \vv^{n}\tau_{\varepsilon}\dd x = \sum_{i,j=1}^d \int_{\Omega}\partial_{x_i\, x_i} b^n \partial_{x_j} b^n \vv^n_j \tau_{\varepsilon} \dd x\\
&=- \sum_{i,j=1}^d \int_{\Omega} \partial_{x_i} b^n \partial_{x_j\, x_i} b^n \vv^n_j \tau_{\varepsilon}+\partial_{x_i} b^n \partial_{x_j} b^n \partial_{x_i} \vv^n_j \tau_{\varepsilon}+\partial_{x_i} b^n \partial_{x_j} b^n \vv^n_j \partial_{x_i} \tau_{\varepsilon} \dd x\\
&=-\int_{\Omega} \frac12 \nabla |\nabla b^n|^2 \cdot \vv^n \tau_{\varepsilon} + (\nabla b^n\otimes \nabla b^n)\cdot \nabla \vv^n \tau_{\varepsilon} +(\nabla b^n\otimes \nabla b^n)\cdot (\vv^n \otimes \nabla \tau_{\varepsilon})\dd x\\
&=-\int_{\Omega} -\frac12 |\nabla b^n|^2  \vv^n \cdot \nabla \tau_{\varepsilon} + (\nabla b^n\otimes \nabla b^n)\cdot \nabla \vv^n \tau_{\varepsilon} +(\nabla b^n\otimes \nabla b^n)\cdot (\vv^n \otimes \nabla \tau_{\varepsilon})\dd x\\
&=-\int_{\Omega} (\nabla b^n\otimes \nabla b^n)\cdot \nabla \vv^n \tau_{\varepsilon}\dd x\\
 &\qquad+\int_{0<\tau_{\varepsilon}<1} \frac12 |\nabla b^n|^2  \vv^n \cdot \nabla \tau_{\varepsilon} -(\nabla b^n\otimes \nabla b^n)\cdot (\vv^n \otimes \nabla \tau_{\varepsilon})\dd x.
\end{aligned}
$$
Thus, thanks to the fact that $b^n \in W^{1,2}(\Omega)$ for almost all time and because of \eqref{tau}, we get
$$
\begin{aligned}
&\int_{\Omega}\Delta b^{n}\,\nabla b^{n}\cdot \vv^{n}\dd x=\lim_{\varepsilon \to 0+} \int_{\Omega}\Delta b^{n}\,\nabla b^{n}\cdot \vv^{n}\tau_{\varepsilon}\dd x
=-\int_{\Omega} (\nabla b^n\otimes \nabla b^n)\cdot \nabla \vv^n \dd x.
\end{aligned}
$$
Hence, returning to \eqref{eq.Gn2maes} and using also \eqref{zero1}, we obtain
\begin{equation}
\int_{\Omega}\frac{|\nabla b^{n}\cdot \vv^{n}|^2}{(b^{n})^2} + \frac{\partial_t b^{n}\, \nabla b^{n}\cdot \vv^{n} }{(b^{n})^2}\dd x = -\int_{\Omega} (\nabla b^n\otimes \nabla b^n)\cdot \nabla \vv^n\dd x.\label{eq.Gn2maes2}
\end{equation}
Finally, we integrate \eqref{eq.Gn2maes2} with respect to time over $(0,t)$ and add the result to \eqref{eq.ID2int} to deduce, after simple algebraic manipulations, that
\begin{equation}\label{eq.ID2int2}
\begin{aligned}
& \frac{\|\nabla b^{n}(t)\|^2_2}{2}+ \int_0^t\int_{\Omega} \frac{|\partial_t b^{n}+\nabla b^{n}\cdot \vv^{n}|^2}{(b^{n})^2}\dd x \dd s\\
&\; \le -\int_0^t\int_{\Omega}(\nabla b^n\otimes \nabla b^n)\cdot \DD^n \dd x\dd s+\int_{\Omega} \ln b^n(t)-\ln b_0 -b^n(t)+b_0\dd x+ {\color{black} \frac{\|b_0\|_{1,2}^2}{2}}.
\end{aligned}
\end{equation}

Now, our goal is to combine this estimate with the balance of the kinetic energy. Therefore, we multiply the $i$-th equation in \eqref{eq.G1n} by $\alpha^n_i(t)$ and sum over $i=1,\ldots, n$ to get
\begin{align*}
\int_{\Omega}\partial_t \vv^{n} \cdot \vv^n - \vv^{n} \otimes \vv^{n} \cdot \nabla \vv^n + \bS^{n}\cdot \nabla \vv^n \dd x&=0.
\end{align*}
Then, for the second term we use integration by parts, so it vanishes thanks to $\divx \vv^n=0$ and for the last term we use \eqref{eq.G3n} to obtain
\begin{align}
\frac12 \frac{\dd}{\dd t} \|\vv^n\|_2^2 +2\|\DD^n\|_2^2 -\int_{\Omega}(\nabla b^n\otimes \nabla b^n)\cdot \DD^n \dd x\dd s&=0.\label{eq.bkn}
\end{align}
Finally, we integrate \eqref{eq.bkn} with respect to time over $(0,t)$ and add the result to \eqref{eq.ID2int2} and obtain the final inequality
\begin{equation}\label{eq.ID2int2*}
\begin{aligned}
& \frac12 \left(\|\vv^n(t)\|_2^2+\|\nabla b^{n}(t)\|^2_2\right)+ \int_0^t\int_{\Omega} \frac{|\partial_t b^{n}+\nabla b^{n}\cdot \vv^{n}|^2}{(b^{n})^2}+2|\DD^n|^2 \dd x \dd s\\
&\qquad \le \int_{\Omega} \ln b^n(t)-\ln b_0 -b^n(t)+b_0\dd x+\frac12 \left( {\color{black}\|b_0\|_{1,2}^2} + \|\vec{P}^n(\vv_0)\|_2^2\right).
\end{aligned}
\end{equation}

To conclude this subsection, we summarize the estimates we have obtained. It directly follows from \eqref{minp}, \eqref{maxp} and \eqref{eq.ID2int2} and Korn's inequality that
\begin{equation}\label{est.uni}
\begin{split}
0<b_{\min}\le b^n \le b_{\max}&<\infty,\\
\sup_{t\in (0,T)} \left(\|\vv^n(t)\|_2 + \|\nabla b^n(t)\|_2\right) &\le C,\\
\int_0^T\|\vv^n\|_{1,2}^2 + \|\partial_t b^{n}+\nabla b^{n}\cdot \vv^{n}\|_2^2\dd t &\le C.
\end{split}
\end{equation}
Consequently, using the interpolation inequality
$$
\|\vv^n\|_{\frac{2(d+2)}{d}}^{\frac{2(d+2)}{d}} \le C \|\vv^n\|^{\frac{4}{d}}_2 \|\vv^n\|^2_{1,2}
$$
we also get from \eqref{est.uni} that
\begin{equation}\label{est.int}
\int_0^T\|\vv^n\|_{\frac{2(d+2)}{d}}^{\frac{2(d+2)}{d}}\dd t\le C.
\end{equation}
Finally, we focus on bounds on time derivatives. First, using \eqref{eq.G1n} and the \textit{a~priori} estimates \eqref{est.uni} we obtain the bound (here we also use the Sobolev embedding $W^{d+1,2}(\Omega)\hookrightarrow W^{1,\infty}(\Omega)$)
\begin{equation}\label{est.tder1}
\int_0^T \|\partial_t \vv^n\|^{\frac{d+2}{d}}_{(W^{1,2}_{0,\divx}\cap W^{d+1,2}(\Omega)^d)^*} \dd t \le C.
\end{equation}
Similarly, with the help of \eqref{est.uni}, \eqref{est.int} and H\"{o}lder's inequality, we deduce that
\begin{equation}\label{est.tder2}
\begin{split}
&\int_0^T \|\partial_t b^n\|_{\frac{d+2}{d+1}}^{\frac{d+2}{d+1}} \dd t\le C\int_0^T \|\partial_t b^n+\nabla b^n \cdot \vv^n\|_{\frac{d+2}{d+1}}^{\frac{d+2}{d+1}} + \||\nabla b^n| |\vv^n|\|_{\frac{d+2}{d+1}}^{\frac{d+2}{d+1}} \dd t\\
&\quad\le C+C\int_0^T \|\partial_t b^n+\nabla b^n \cdot \vv^n\|_{2}^{2} + \|\nabla b^n\|_2^2 + \|\vv^n\|_{\frac{2(d+2)}{d}}^{\frac{2(d+2)}{d}} \dd t\le C.
\end{split}
\end{equation}

\subsection{The limit $n\to \infty$}
In this final subsection, we let $n\to \infty$ and show that the sequence of solutions to the approximating problem indexed by $n$ tends to a solution of the original problem. Having the uniform estimates \eqref{maxp}, \eqref{minp}, \eqref{est.uni} and \eqref{est.int}--\eqref{est.tder2}, there exist $\vv$ and $b$ and subsequences that we do not relabel, such that
\begin{align}
\vv^n &\rightharpoonup^* \vv &&\textrm{weakly$^*$ in } L^{\infty}(0,T; L^2_{0,\divx}), \label{n1}\\
\vv^n &\rightharpoonup \vv &&\textrm{weakly in } L^{2}(0,T; W^{1,2}_{0,\divx})\cap L^{\frac{2(d+2)}{d}}(0,T; L^{\frac{2(d+2)}{d}}(\Omega)^d), \label{n2}\\
\partial_t\vv^n &\rightharpoonup \partial_t\vv &&\textrm{weakly in } L^{\frac{d+2}{d}}(0,T;(W^{1,2}_{0,\divx}\cap W^{d+1,2}(\Omega)^d)^*), \label{n3}\\
b^n &\rightharpoonup^* b &&\textrm{weakly$^*$ in } L^{\infty}(0,T; W^{1,2}(\Omega)), \label{n4}\\
b^n &\rightharpoonup^* b &&\textrm{weakly$^*$ in } L^{\infty}(0,T; L^{\infty}(\Omega)), \label{n5}\\
\partial_t b^n &\rightharpoonup \partial_t b &&\textrm{weakly in } L^{\frac{d+2}{d+1}}(0,T; L^{\frac{d+2}{d+1}}(\Omega)). \label{n6}
\end{align}
Using the Aubin--Lions lemma we then deduce that
\begin{align}
\vv^n &\to \vv &&\textrm{strongly in } L^{p}(0,T; L^{p}(\Omega)^d)\label{n7}\\
\intertext{for all  $p\in [1,2(d+2)/d)$. {\color{black} In addition, }}
\partial_t b^n + \nabla b^n \cdot \vv^n &\rightharpoonup \partial_t b +\nabla b \cdot \vv &&\textrm{weakly in } L^{2}(0,T; L^{2}(\Omega)), \label{n8}\\
b^n &\to b &&\textrm{strongly in } \mathcal{C}([0,T]; L^q(\Omega))\label{n9}
\end{align}
for all $q\in [1,\infty)$. {\color{black} In particular, \eqref{n9} follows from \eqref{n4} and \eqref{n6}, the compact embedding of $L^{\infty}(0,T; W^{1,2}(\Omega)) \cap W^{1,\frac{d+2}{d+1}}(0,T;L^{\frac{d+2}{d+1}}(\Omega))$ into $\mathcal{C}([0,T]; L^q(\Omega))$ for all $q \in [1,\infty)$ when $d=2$ and all $q \in [1,6)$ when $d=3$, guaranteed by the Aubin--Lions lemma; and, in the case of $d=3$, we have also used \eqref{n5} in order to extend the range of $q$ from $[1,6)$ to $[1,\infty)$ using the boundedness of the sequence
$(b^n)_{n \geq 1}$ in $L^\infty(0,T;L^\infty(\Omega))$, guaranteed by \eqref{n5}, in conjunction with the strong convergence
of this sequence in $\mathcal{C}([0,T]; L^q(\Omega))$ for $q \in [1,6)$.}

In addition, it follows from \eqref{minp}, \eqref{maxp} and \eqref{n9} that almost everywhere in $Q$ we have
\begin{equation}
0<b_{\min}\le b(t,x) \le b_{\max} < \infty. \label{minimax}
\end{equation}
Furthermore, we can easily let $n\to \infty$ in \eqref{eq.Gn2m} to obtain, for almost all $t\in(0,T)$, that
\begin{equation}
\int_{\Omega} \frac{\partial_t b\, w}{b^2} + \frac{\nabla b \cdot \vv}{b^2}w + (1 - b^{-1})w + \nabla b\cdot \nabla w \dd x  = 0 \qquad \textrm{ for all } w\in W^{1,2}(\Omega),\label{eq.Gn2fin}
\end{equation}
which is nothing else but \eqref{Tw2} with our special choices of $\nu_1$, $\mu$ and $\sigma$.
To obtain also the strong convergence of $\nabla b^n$, we set $w:= b$ in \eqref{eq.Gn2fin} and $w:=b^n$ in \eqref{eq.Gn2m} and with the help of \eqref{n8} and \eqref{n9}, we deduce that
$$
\begin{aligned}
\lim_{n \to \infty}\int_0^T \|\nabla b^n\|_2^2 \dd t &=-\lim_{n\to \infty}\int_{Q} \frac{\partial_t b^n +\nabla b^n \cdot \vv^n}{b^n} + (b^n - 1) \dd x \dd t\\
&=-\int_{Q} \frac{\partial_t b +\nabla b \cdot \vv}{b} + (b - 1) \dd x \dd t =\int_0^T \|\nabla b\|_2^2\dd t.
\end{aligned}
$$
Thus, this identity, combined also with \eqref{n4}, leads to
\begin{align}
b^n &\to b &&\textrm{strongly in } L^{2}(0,T; W^{1,2}(\Omega)). \label{n10}
\end{align}
Consequently, by \eqref{n1}--\eqref{n3}, \eqref{n7} and \eqref{n10}, we can let $n\to \infty$ in \eqref{eq.G1n}, \eqref{eq.G3n} to conclude that for almost all $t\in (0,T)$ and all $\vw \in W^{1,2}_{0,\divx} \cap W^{d+1,2}(\Omega)^d$
we have
\begin{align}
\langle \partial_t \vv,  \vw\rangle  +\int_{\Omega}  \bS\cdot \nabla \vw - (\vv \otimes \vv) : \nabla \vw\dd x&=0,\label{eq.G1fin}
\end{align}
with $\bS$ given as
\begin{equation}
\bS=  2\DD - (\nablax b \otimes \nablax b), \quad \DD=\frac{1}{2}(\nabla \vv+(\nabla \vv)^{\textrm{T}}). \label{eq.G3fin}
\end{equation}
Hence, \eqref{Tw1} is also satisfied. Furthermore, letting $n\to \infty$ in \eqref{eq.ID2int2}, using lower semicontinuity of norms and the convergence results \eqref{n1}, \eqref{n2}, \eqref{n4}, \eqref{n8} and \eqref{n9}, we also obtain the energy inequality \eqref{T:EI}. Thus, to complete the proof of Theorem~\ref{T:main}  it only remains to prove the attainment of the initial conditions.

First, a standard argument yields that
\begin{equation}\label{weakatt}
\begin{aligned}
\vv(t)&\rightharpoonup \vv_0 &&\textrm{ weakly in } L^2_{0,\divx},\\
b(t) &\rightharpoonup b_0 &&\textrm{ weakly in }L^2(\Omega).
\end{aligned}
\end{equation}
Next, thanks to \eqref{n4} and \eqref{n6}, we can strengthen the above convergence results as follows:
\begin{equation}\label{weakatt2}
\begin{aligned}
\nabla b(t)&\rightharpoonup \nabla b_0 &&\textrm{ weakly in } L^2(\Omega)^d,\\
b(t) &\to b_0 &&\textrm{ strongly in }L^2(\Omega).
\end{aligned}
\end{equation}
In addition, using \eqref{n9} and weak lower semicontinuity of the norm function and neglecting the terms with the correct sign, we can let $n\to \infty$ in \eqref{eq.ID2int2} to obtain for all $t\in (0,T)$ that
\begin{equation}\label{eq.ID2intfin}
\begin{aligned}
& \frac{\|\vv(t)\|_2^2+\|\nabla b(t)\|^2_2}{2} \le \int_{\Omega} \ln \left(\frac{b(t)}{b_0}\right) -b(t)+b_0\dd x+ \frac{\|\nabla b_0\|_2^2+ \|\vv_0\|_2^2}{2}.
\end{aligned}
\end{equation}
Hence, it follows (by noting the strong convergence results \eqref{weakatt2} and \eqref{minimax}) that
$$
\limsup_{t\to 0}  \|\vv(t)\|_2^2+\|\nabla b(t)\|^2_2 \le  \|\vv_0\|_2^2+\|\nabla b_0\|^2_2,
$$
which, when combined with \eqref{weakatt} and \eqref{weakatt2}, then leads to \eqref{T:ic}. The proof of the theorem is thereby complete.

\begin{appendix}
\section{Existence of Galerkin approximation}
Here, we sketch the proof of the existence of the Galerkin approximation, i.e., a solution to \eqref{eq.G1}--\eqref{eq.G3}. First, we recall the problem we would like to investigate, but we omit writing all superscripts related to the indices of the approximating sequences. Hence, we are interested in finding a solution to
\begin{align}
\begin{aligned}
\int_{\Omega}\partial_t \vv \cdot \vw_i - \vv \otimes \vv : \nabla \vw_i + \bS : \nabla \vw_i &\dd x = 0, \qquad i=1,\ldots, n,
\end{aligned} \label{eq.A1} \\
\begin{aligned}
\int_{\Omega} \frac{\partial_t b\, w_j}{(T_n(b))^2} &+ \frac{\nabla b\cdot \vv}{(T_n(b))^2}w_j
+ (1 - (T_n(b))^{-1})w_j \dd x \\
 &+ \int_{\Omega}  \nabla b\cdot \nabla w_j \dd x  = 0, \qquad j=1,\ldots, \ell, \end{aligned} \label{eq.A2}
\end{align}
where $\bS^{n,\ell}$ is defined a.e. in $Q$ by
\begin{equation}
\bS=  2\DD - (\nablax b \otimes \nablax b), \quad \DD=\frac{1}{2}\big(\nabla \vv+(\nabla \vv)^{\textrm{T}}\big), \label{eq.A3}
\end{equation}
where $(\vv,b)$ are given by
$$
\vv(t,x):= \sum_{i=1}^n \alpha_{i}(t) \vw_i(x),\qquad b(t,x):= \sum_{i=1}^{\ell} \beta_i(t) w_i(x)
$$
with the initial conditions
$$
\vv(0,x)= \vec{P}^n(\vv_0)(x), \qquad b(0,x)= P^{\ell}(b_0)(x).
$$

Because of the presence of the nonlinearity in the term involving $\partial_t b$ we cannot simply invoke standard existence results from the theory of ordinary differential equations; we therefore proceed slightly differently. We denote
$$
f(s):=\int_0^s \frac{1}{(T_n(t))^2} \dd t, \qquad F(s):= \int_{0}^{s} f(t)\dd t
$$
and observe that $F$ is uniformly convex, satisfying, with some $C_1,C_2>0$,
\begin{equation}
C_1 s^2 \le F(s) \le C_2 s^2, \qquad s \geq 0.\label{FF}
\end{equation}
In addition, we can rewrite \eqref{eq.A2} as
\begin{align}
\begin{aligned}
\int_{\Omega} \partial_t f(b)\,w_j &+ \nabla f(b)\cdot \vv w_j
+ (1 - (T_n(b))^{-1})w_j \dd x \\
 &+ \int_{\Omega}  \nabla b\cdot \nabla w_j \dd x  = 0, \qquad j=1,\ldots, \ell. \end{aligned} \label{eq.B2}
\end{align}
We shall approximate \eqref{eq.B2} with the help of Rothe's method, i.e., we consider the following system of $N \times \ell$ equations:
\begin{align}
\begin{aligned}
&\int_{\Omega} \frac{f(b^{k+1})-f(b^k)}{\tau} w_j + \nabla f(b^{k+1})\cdot \vv^k w_j
+ (1 - (T_n(b^{k+1}))^{-1})w_j \dd x \\
 &+ \int_{\Omega}  \nabla b^{k+1}\cdot \nabla w_j \dd x  = 0, \qquad j=1,\ldots, \ell, \qquad k=0,\ldots, N-1, \end{aligned} \label{eq.BB2}
\end{align}
where $\tau:=T/N$ and
$$
b^0(x):= P^{\ell}b_0 (x), \quad \vv^k(x):=\frac{1}{\tau} \int_{t_k}^{t_{k+1}} \vv(t,x) \dd t, \quad t_k:=kT/N,
$$
and
$$
b^k(x):= \sum_{i=1}^{\ell} \beta^k_i w_i(x).
$$
Finally, we define
\begin{equation}\label{BBB}
b(t,x):=\frac{t_{k+1}-t}{\tau} b^k(x) + \frac{t-t_{k}}{\tau}b^{k+1}(x) \qquad \textrm{ for } t\in [t_k, t_{k+1}].
\end{equation}

Notice that for given $\vv$, the problem \eqref{eq.BB2} has {\color{black} a solution}.\footnote{\color{black} The proof of existence of a solution $b^{k+1} \in W_\ell:=\mbox{span}\{w_i\}_{i=1}^\ell$ proceeds, by using a simple consequence of Brouwer's fixed point theorem, as follows. For a given, fixed, $\vv^k$ as above and a given, fixed, $b^k$, with $k \in \{0,\dots, N-1\}$,
we define the (continuous) mapping  $\Phi_k\,:\, W_\ell \rightarrow W_\ell$ by
\begin{align*}
\begin{aligned}
&(\Phi_k(b),w):=\int_{\Omega} \frac{f(b)-f(b^k)}{\tau} w + \nabla f(b)\cdot \vv^k w
+ (1 - (T_n(b))^{-1})w \dd x
+ \int_{\Omega}  \nabla b \cdot \nabla w \dd x, \end{aligned}
\end{align*}
for all $w \in W_\ell$. An argument analogous to the one leading to \eqref{eq.BB2a} results in
\begin{align*}
\tau (\Phi_k(b),b)=\int_{\Omega}(f(b)-f(b^k))\,b + \tau(1 - (T_n(b))^{-1})\,b\dd x + \tau\|\nabla b\|_2^2\qquad \forall\, b \in W_\ell.
\end{align*}
As $(1 - (T_n(s))^{-1})\,s \geq -1$ for all $s \in \mathbb{R}$, it follows that $\tau (\Phi_k(b),b) \geq (f(b)-f(b^k), b) - \tau |\Omega|$. By noting that $T_n(s) \leq n$ for all $s \in \mathbb{R}$ and recalling the definition of $f$, we have $$\tau (\Phi_k(b),b) \geq (1/n^2)\|b\|^2_2 - (f(b^k),b) - \tau |\Omega|
\geq (1/2n^2) \|b\|^2_2 - (n^2/2)\|f(b^k)\|^2_2 - \tau |\Omega|.$$ Hence, $(\Phi_k(b),b)\geq 0$ for all $b \in W_\ell$ such that
$\|b\|^2 \geq  n^4 \|f(b^k)\|^2_2 + 2\tau n^2 |\Omega|$, and the existence of a $b=b^{k+1} \in W_\ell$ such that $\Phi_k(b)=0$
then follows from  Corollary 1.1 on p.279 in \cite{Gir-Rav}.
}
Moreover, this solution and consequently also $b$ depends continuously on $\vv$ and therefore, we can now use Carath\'{e}odory's theory to deduce the existence of a solution to the problem {\color{black}\eqref{eq.A1}--\eqref{eq.A3}} and \eqref{eq.BB2}, \eqref{BBB}. Finally, we will establish \textit{a priori} estimates in which we can easily let $N \to \infty$ to deduce the solvability of the original problem.

Multiplying \eqref{eq.BB2} by $\beta^{k+1}_j$, summing with respect to $j$ and using $\divx \vv =0$, we obtain
\begin{align}
\begin{aligned}
&\int_{\Omega}(f(b^{k+1})-f(b^k))\,b^{k+1} + \tau(1 - (T_n(b^{k+1}))^{-1})\,b^{k+1} \dd x + \tau\|\nabla b^{k+1}\|_2^2 = 0 \end{aligned} \label{eq.BB2a}
\end{align}
for all $k=0,\ldots, N-1$. Next, we denote by $F^*$ the convex conjugate of the uniformly convex function $F$, i.e.,
$$
F^*(s):= \sup_{t} \left(st - F(t)\right), \qquad \textrm{ with } \tilde{C}_1s^2\le F^*(s)\le \tilde{C}_2 s^2,
$$
and with the help of Young's inequality we have
$$
\begin{aligned}
(f(b^{k+1})-f(b^k))b^{k+1}&=F(b^{k+1}) + F^*(f(b^{k+1})) - f(b^k) b^{k+1}\\
&\ge F^*(f(b^{k+1})) - F^*(f(b^{k})).
\end{aligned}
$$
{\color{black} We note further that, since $1/n \leq T_n(s) \leq n$ for all $s \in \mathbb{R}$, we have $f(s) \geq s/n^2$ for all $s \geq 0$; note furthermore that if $b^{k+1}(x) \leq 0$ for some $x \in \Omega$ then $(1 - (T_n(b^{k+1}(x)))^{-1})\,b^{k+1}(x) \geq 0$. Hence, by using Young's inequality and the properties of the function $F^*$, we have that
\begin{align*}
&\int_\Omega \tau(1 - (T_n(b^{k+1}))^{-1})\,b^{k+1} \dd x \geq \int_\Omega \tau(1 - (T_n(b^{k+1}))^{-1})\,b^{k+1}
\chi_{b^{k+1} \geq 0} \dd x
\\
&\qquad\geq -  \frac{\tau n^4}{2\tilde{C}_1}\int_\Omega (1 - (T_n(b^{k+1}))^{-1})^2 \chi_{b^{k+1} \geq 0} \dd x- \frac{\tilde{C}_1\tau }{2n^4} \int_\Omega (b^{k+1})^2\chi_{b^{k+1} \geq 0}  \dd x\\
&\qquad \geq - \frac{\tau n^4}{2\tilde{C}_1}(n-1)^2 |\Omega| - \frac{\tilde{C}_1\tau }{2} \int_\Omega (f(b^{k+1}))^2 \chi_{b^{k+1} \geq 0} \dd x\\
&\qquad \geq - \frac{\tau n^4}{2\tilde{C}_1}(n-1)^2 - \frac{\tau}{2} \int_\Omega F^*(f(b^{k+1}))\, \chi_{b^{k+1} \geq 0} \dd x\\
&\qquad \geq - \frac{\tau n^4}{2\tilde{C}_1}(n-1)^2 - \frac{\tau}{2} \int_\Omega F^*(f(b^{k+1})) \dd x.
\end{align*}
}

\noindent
Thus, using these inequalities in \eqref{eq.BB2a} and summing the result over $k=1, \ldots, m$ with arbitrary $m\le N-1$, {\color{black} we get,  with $C=T n^4 (n-1)^2/(2\tilde{C}_1)$, that}
\begin{align}
\begin{aligned}
&\int_{\Omega}F^*(f(b^{m+1}))- F^*(f(b^0)) \dd x + \sum_{k=0}^m \tau \|\nabla b^{k+1}\|_2^2\\
&\qquad \le C + \frac12\sum_{k=0}^m \tau \int_{\Omega} F^*(f(b^{k+1}))\dd x.\end{aligned} \label{eq.BB2b}
\end{align}
Hence, using the discrete version of the Gronwall's lemma and also the definition of $b$, we obtain
\begin{equation}\label{A.10}
\sup_{t\in (0,T)} \|b(t)\|_2^2 + \int_0^T \|\nabla b(t)\|_2^2 \dd t \le C.
\end{equation}
With this estimate in hand, we can now also obtain an estimate for $\vv$ (see the proof of the main theorem):
\begin{equation}\label{A.11}
\sup_{t\in (0,T)} \|\vv(t)\|_2^2  \le C(n,\ell) \implies \sup_{t\in (0,T)} \|\nabla \vv(t)\|_{\infty}^2\le C(n,\ell).
\end{equation}
Finally, we multiply \eqref{eq.BB2} by $(\beta^{k+1}_j-\beta^k_j)$ and  sum with respect to $j$ to obtain
\begin{align}
\begin{aligned}
&\int_{\Omega}\frac{(f(b^{k+1})-f(b^k))(b^{k+1}-b^k)}{\tau}  + (1 - (T_n(b^{k+1}))^{-1})(b^{k+1}-b^k) \dd x \\
&\qquad +\int_{\Omega} \nabla b^{k+1} \cdot \nabla (b^{k+1}-b^k)\dd x = 0. \end{aligned} \label{eq.BB2d}
\end{align}
Thus, summing with respect to $k=0,\ldots, m$ and using standard inequalities, we have
\begin{align}
\begin{aligned}
&\sum_{k=0}^m \tau \left\| \frac{b^{k+1}-b^k}{\tau}\right\|_2^2  +\|\nabla b^m\|_2^2 \le C(1+\|\nabla b^0\|_2^2). \end{aligned} \label{eq.BB2e}
\end{align}
Consequently, using the definition of $b$ we get
\begin{equation}
\int_0^T \|\partial_t b\|_2^2 \dd t + \sup_{t\in (0,T)} \|\nabla b(t)\|_2^2 \le C.\label{konec}
\end{equation}
Therefore, with the help of \eqref{A.10}, \eqref{A.11} and \eqref{konec}, it is quite easy to let $N\to \infty$ to establish the existence of a solution to the original problem.

\end{appendix}



\bibliographystyle{amsplain}
\bibliography{vit-prusa}

\end{document}